\documentclass[11pt,reqno]{amsart}

\usepackage{amsmath, amsfonts, amssymb, amscd, enumerate, color, graphics}
\theoremstyle{plain}
\newtheorem{theo}{Theorem}[section]
\newtheorem{lem}[theo]{Lemma}

\theoremstyle{definition}
\newtheorem{rem}[theo]{Remark}

\newtheorem{definition}[theo]{Definition}

\newenvironment{pf}{\noindent{\it Proof. }}{$\square$\par\medskip}
\newenvironment{pfns}{\noindent{\it Proof. }}{\par\medskip}

\newenvironment{colored}{\color{red}}{}
\newcommand{\bc}{\begin{colored}}
\newcommand{\ec}{\end{colored}}

\theoremstyle{plain}
\newtheorem{lemma}[theo]{Lemma}

\theoremstyle{definition}

\renewcommand{\=}{\overset{\operatorname{def}}{=}}

\newcommand{\beq}{\begin{equation}}
\newcommand{\eeq}{\end{equation}}
\renewcommand{\a}{\alpha}
\renewcommand{\b}{\beta}

\newcommand{\g}{\gamma}
\newcommand{\h}{\eta}

\renewcommand{\l}{\lambda}
\newcommand{\m}{\mu}
\renewcommand{\o}{\omega}

\renewcommand{\r}{\rho}
\newcommand{\s}{\sigma}
\renewcommand{\t}{\tau}

\renewcommand{\L}{\Lambda}
\renewcommand{\O}{\Omega}

\newcommand{\bC}{\mathbb{C}}

\newcommand{\bR}{\mathbb{R}}

\renewcommand{\gg}{\mathfrak{g}}
\newcommand{\gh}{\mathfrak{h}}

\newcommand{\gl}{\mathfrak{l}}
\newcommand{\gm}{\mathfrak{m}}
\newcommand{\gn}{\mathfrak{n}}

\newcommand{\so}{\mathfrak{so}}

\newcommand{\ggl}{\mathfrak{gl}}
\newcommand{\ggr}{\mathfrak{gr}}

\newcommand\GL{\mathrm{GL}}

\newcommand\SO{\mathrm{SO}}

\newcommand{\cB}{\mathcal{B}}
\newcommand{\cC}{\mathcal{C}}
\newcommand{\cD}{\mathcal{D}}
\newcommand{\cE}{\mathcal{E}}
\newcommand{\cF}{\mathcal{F}}

\newcommand{\cH}{\mathcal{H}}
\newcommand{\cI}{\mathcal{I}}

\newcommand{\cK}{\mathcal{K}}
\newcommand{\cL}{\mathcal{L}}

\newcommand{\cT}{\mathcal{T}}
\newcommand{\cU}{\mathcal{U}}
\newcommand{\cV}{\mathcal{V}}

\renewcommand{\square}{\kern1pt\vbox
{\hrule height 0.6pt\hbox{\vrule width 0.6pt\hskip 3pt
\vbox{\vskip 6pt}\hskip 3pt\vrule width 0.6pt}\hrule height0.6pt}\kern1pt}

\DeclareMathOperator\Aut{Aut}

\DeclareMathOperator\Ad{Ad}
\DeclareMathOperator\ad{ad}

\DeclareMathOperator\Id{Id}

\DeclareMathOperator{\Tor}{Tor}

\renewcommand\Re{\operatorname{Re}}
\renewcommand\Im{\operatorname{Im}}

\newcommand{\Hom}{{\operatorname{Hom}}}

\newcommand{\wt}{\widetilde}
\newcommand{\wh}{\widehat}

\newcommand{\be}{\begin{equation}}
\newcommand{\ee}{\end{equation}}

\def\<#1,#2>{\langle\,#1,\,#2\,\rangle}
\newcommand{\arr}{\begin{array}{rlll}}
\newcommand{\ea}{\end{array}}
\newcommand{\bea}{\begin{eqnarray}}
\newcommand{\eea}{\end{eqnarray}}
\newcommand{\bean}{\begin{eqnarray*}}
\newcommand{\eean}{\end{eqnarray*}}

\def\sideremark#1{\ifvmode\leavevmode\fi\vadjust{
\vbox to0pt{\hbox to 0pt{\hskip\hsize\hskip1em
\vbox{\hsize3cm\tiny\raggedright\pretolerance10000
\noindent #1\hfill}\hss}\vbox to8pt{\vfil}\vss}}}

\newcounter{ssig}
\newcounter{ttig}
\newcommand{\under}[1]{{\underline{#1}\,}}
\renewcommand{\Vert}{\operatorname{Vert}}

\title[The equivalence problem for Levi degenerate CR manifolds]
{The equivalence problem for 5-dimensional  \\Levi degenerate CR manifolds}

\address{
Costantino Medori,
Dipartimento di Matematica e Informatica,
Universit\`a di Parma,
Parma, Italy.}
\email{costantino.medori@unipr.it}
\address{
Andrea Spiro, Scuola di Scienze e Tecnologie, Universit\`a di Camerino, 
Camerino, 
Italy.}
\email{andrea.spiro@unicam.it}


\author{Costantino Medori and Andrea Spiro}

\begin{document}
\begin{abstract} Let $M$ be a CR manifold of hypersurface type, which is Levi degenerate but  also satisfying a $k$-nondegeneracy condition  at all points.  This  might be only if $\dim M \geq 5$ and  if  $\dim M = 5$, then   $k= 2$ at all points.  We prove that for any  5-dimensional, uniformly 2-nondegenerate CR manifold $M$  there exists  a canonical   Cartan connection, modelled on a suitable projective completion  of the tube over the future light cone  $\{z \in \bC^3: (x^1)^2 + (x^2)^2 - (x^3)^2 = 0\,,\,x^3 > 0\}$. 
This   determines a complete  solution to the equivalence problem for this class of  CR manifolds. 
\end{abstract}
\maketitle
\null \vspace*{-.25in}
\section{Introduction}
Let $M$ be a 5-dimensional CR hypersurface,  which is Levi degenerate at all points.  Quite simple examples are  provided  by  Cartesian products of the form    $M = \overline M \times \bC$ for some 3-dimensional CR manifold $\overline M$. A much less trivial case  is represented  by the so-called {\it tube over the future light cone}  
\beq \label{tube0} \cT = \{z \in \bC^3: (x^1)^2 + (x^2)^2 - (x^3)^2 = 0\,,\,x^3 > 0\} \subset \bC^3\,.\eeq
This  hypersurface  is  in fact Levi degenerate at all points (it is homogeneous w.r.t. Ê$\Aut(\cT)$), it is foliated by complex leaves  and yet it admits no local  {\it CR straightening}, that is no local   CR equivalence with a Cartesian product  of the form $\overline M \times \bC$.\par
\smallskip
Freeman (\cite{Fr}) found  necessary and sufficient conditions for real analytic CR manifolds to admit local straightenings, together with  obstructions to the existence of CR straightenings in the   smooth category.  Such obstructions  are equivalent  to the so-called {\it $k$-nondegeneracy conditions}  at  its points (\cite{BER, KZ}).  We recall that a CR hypersurface $M$  satisfies the $1$-nondegeneracy condition at all points  if and only if it is Levi nondegenerate and    that   the   other $k$-nondegeneracy conditions  for  $k \geq 2$ can be taken as progressively weaker nondegeneracy conditions. 
\par
\smallskip
The smallest possible dimension for a CR hypersurface $M$ to be  Levi degenerate and yet $k$-nondegenerate at all points  is $5$. In  such a case $k$ is necessarily   equal to $2$.  For brevity, we  call the 5-dimensional, 2-nondegenerate CR hypersurfaces  of uniform type  {\it girdled CR manifolds}.\par
\smallskip
These CR manifolds have been recently considered in several studies, which for instance brought to    Ebenfelt's  normal forms for   real analytic, Levi degenerate  CR hypersurfaces  in $\bC^3$ and Kaup and Fels'   classification  of homogeneous Levi degenerate 5-dimensional CR manifolds  (\cite{Eb, FK, FK1}).\par
\smallskip
In this paper, we give a general solution to the equivalence problem for girdled CR manifolds in the $\cC^\infty$ category,  proving the existence of a canonical Cartan connection for  any such manifold $M$.  \par
\smallskip
We recall that (see e.g. \cite{Kb, Sh})  a  {\it  Cartan connection  on a manifold $N$, modelled on a homogeneous space $G/H$,}  is 
a pair $(Q, \varpi)$,  formed by 
 a principal $H$-bundle $\pi: Q \longrightarrow N$  and a $\gg$-valued 1-form 
 $\varpi:TQ \longrightarrow \gg = Lie(G)$ 
satisfying the following conditions:
\begin{itemize}
\item[(a)]  $\varpi_y: T_y Q \longrightarrow  \gg$ is a linear isomorphisms for any $y\in Q$  and
$$(\varpi_y)^{-1}|_{\gh}: \gh = Lie(H) \longrightarrow T^{\Vert}_y Q$$
is the standard  isomorphism given by the  action of $H$ on $Q$, 
\item[(b)] $R^*_h \varpi = \Ad_{h^{-1}} \varpi$  for any $h\in H$. 
\end{itemize}
If $N$ is a manifold endowed with  a fixed geometric structure,  a Cartan connection  $(Q, \varpi)$  is called  {\it canonical} if there exists a natural  correspondence between 
the automorphisms  of the geometric structure  and the  diffeomorphisms $\wh f: Q \longrightarrow Q$ such that
$\wh f^* \varpi = \varpi$. We point out that, if $(Q, \varpi)$ is a  canonical Cartan connection,  any fixed  basis $(E^o_i)$ of $\gg$  gives a canonical  {\it absolute parallelism} on $Q$ (also called {\it $\{e\}$-structure}),  namely the collection of vector fields $$E_i|_y = \varpi_y^{-1}(E^o_i)\,,\qquad y \in Q\,.$$
Since  the structure functions of a canonical  absolute parallelism  give a complete set of invariants for the  geometric structure on $N$  (see e.g. \cite{St, AS}),   for the geometric structures that  admit canonical   Cartan connections, the equivalence problems have  exact and complete solutions. \par
\smallskip
Canonical Cartan connections  give also  valuable information on the  automorphism groups of  the considered geometric structures and allow constructions  of   useful special coordinates (see e.g.  \cite{Kb, Sh, SS1}).\par
\smallskip
Our main result on girdled CR manifolds is the   following:
\begin{theo} For any  5-dimensional girdled CR manifold  $(M, \cD, J)$,   there exists a canonical 
Cartan connection $(Q, \varpi)$,  modelled on Ê the projective completion $M_o = \SO_{3,2}^o/H \subset \bC P^4$ of   $\cT$ (see  \S \ref{homogeneous model} for definition of $M_o$).  
\end{theo}
The proof is constructive and provides an explicit description of the  bundle $\pi: Q\longrightarrow M$ and of the $\gg$-valued 1-form $\varpi$. Roughly speaking, it  consists of 
a suitable modification  of   Tanaka's construction of  Cartan connections for  geometric structures modelled on semi-simple Lie groups (\cite{Ta3, AS}).  In particular, our result can be taken as    the  analogue of  the   Cartan connections 
of Levi nondegenerate   hypersurfaces (Tanaka's and Chern-Moser's connections) and  of  the CR manifolds endowed  with the so-called {\it parabolic geometries} (\cite{Ta2, Ta3, CM, CS, SSl, SS, SS1}).  \par
Note that, since the  isotropy   of the model space Ê$M_o = \SO_{3,2}^o/H$  is not a parabolic subgroup,  the girdled CR structures cannot be classified as  parabolic geometry. On the other hand,  our modifications  of Tanaka's approach can be very likely extended  to other dimensions and other types of finitely nondegenerate  CR structures and we expect the existence of an  interesting new  class of geometries,  which  includes both girdled CR structures and  parabolic geometries as special cases. \par
\smallskip
We  conclude mentioning  that another  absolute parallelism for  girdled CR manifolds was previously  determined by Ebenfelt   in \cite{Eb1}. However, such parallelism  can be considered  as canonical   only if one  considers  CR automorphisms satisfying certain additional assumptions. It consequently  provides only solutions   to   the corresponding  restricted equivalence problem. \par
\smallskip
After finishing  our paper, we realised  that  Isaev and Zaitsev   recently constructed an  
absolute parallelism for girdled  CR manifolds, in general not determined   by a Cartan connection,  which  therefore gives an alternative  solution to the   equivalence problems for such manifolds (\cite{IZ}). \par 
\medskip
The paper is organised as follows: in \S 2 and \S 3 we give the basic  definitions and  properties  of girdled CR manifolds, of the tube $\cT$ over the future light cone and of its projective completion $M_o \subset \bC P^4$;  in \S 4, we introduce  some definitions and simple facts  on vector spaces with  filtrations, which will be  used in  later sections; in  \S 5, \S 6 and \S 7, we construct the three steps of a tower, which is canonically associated with a girdled CR manifold and is   the analogue of Tanaka's tower of Levi-nondegenerate CR manifolds;  in  \S 8, we prove the main theorem and determine  the structure equations of   a girdled manifold. \par
\medskip
{\it Notation.} In the following, for any given (real or complex)  subbundle $\cK \subset T^\bC M$ of the complexified tangent space $T^\bC M$,  we  indicate by 
$\under \cK$ the class of all local smooth sections of  $\cK$.\par
\smallskip
Ê{\it Acknowledgement.} We are grateful to A. Isaev for pointing out  a mistake in our previous description of  $\Aut(\cT)$ and $\Aut(M_o)$. 
\section{Preliminaries}
\setcounter{equation}{0}
\subsection{Finitely nondegenerate CR manifolds of dimension $5$} \hfill\par
Given a $2n+k$-dimensional manifold $M$, a {\it CR structure  on $M$ of codimension $k$} is  a pair $(\cD, J)$,  formed by  a distribution $\cD \subset TM$ of codimension $k$ and 
a smooth family of complex structures $J_x: \cD_x \longrightarrow \cD_x$ satisfying the following integrability condition: 
\begin{itemize}
\item[]  the  bundle $\cD^{10} \subset T^\bC M$,  given by the $+i$-eigenspaces $\cD^{10}_x \subset \cD^\bC_x$ of  the complex structure   $J_x$,   is involutive, i.e., 
$$[X^{10}, Y^{10}] \in \under \cD^{10}\quad \text{for any pair}\,X^{10}, Y^{10} \in \under \cD^{10}\,.$$
\end{itemize}
The complex vector bundles $\cD^{10}$ and  $\cD^{01} = \overline{\cD^{10}}$   are called {\it holomorphic and anti-holomorphic  bundles} of the  CR structure $(\cD,J)$, respectively.
\par
Two CR manifolds $(M, \cD, J)$ and $(M', \cD', J')$ are called {\it (locally) CR equivalent} if there exists a (local) diffeomorphism $f: M \longrightarrow M'$ such that
 $$f_*(\cD) = \cD'\,,\qquad f_*(J) = J'\,.$$ \par
\smallskip
  The  {\it  Freeman sequence} of  a CR manifold $(M, \cD, J)$ (\cite{Fr}, Thm. 3.1) is the nested sequence of families of  complex vector fields
$$\dots  \subset  \under\cF_{ j+1} \subset  \under\cF_{ j} \subset \ldots \subset   \under \cF_{1} \subset \under \cF_{0} \subset \under \cF_{-1} = \under \cD^\bC$$
iteratively defined by 
\beq \label{1.1}
\begin{array}{l}
\under \cF_{k+1}\!\! =\under \cF_{k+1}^{10} + 
\overline{\under \cF_{k+1}^{10}}\qquad \text{with} 
\ \,  \under{\cF}^{10}_{-1} = \under{\cD}^{10} \ \,\text{and}
\\ \,\\
 \phantom{aaaaaaaaaaa} \under{\cF}^{10}_{k+1}=\{X\in \under{\cF}^{10}_k :[X, \under{\cD}^{01}] = 0\!\!\!\!\! \mod  \under{\cF}^{10}_k + \under{\cD}^{01} \}.
\end{array}
\eeq
 A CR structure $(\cD, J)$  is called  {\it   regular} if  the vector fields  in  $\under{\cF}_j$ and $\under{\cF}^{10}_j$  are the sections of  corresponding  complex distributions  $\cF_j, \cF^{10}_j \subset \cD^\bC$ for any $j \geq -1$.  From now on, {\it any CR manifold will be tacitly assumed to be   regular}.  \par
 \medskip
 The  following  is a  consequence of  definitions.  \par
\begin{lem}\label{collectanea} If $(M, \cD, J)$ is a regular CR manifold, all complex distributions $\cF_k$ are $J$-invariant and real (i.e.  equal to their  conjugate). \par
Furthermore, the class of  vector fields in $\cE = \Re(\cF_0) \subset \cD$  is equal to
\beq \label{definitionE} \under \cE =  \{\,X \in \under \cD\,:\,[X, \under \cD] \subset \under \cD\} \,.\eeq
In particular, $\cE$ is an involutive subdistribution of $\cD$. 
\end{lem}
We may now consider the following definition (see e.g. \cite{BER, KZ}).  \par
\begin{definition}  A (regular) CR manifold $(M, \cD, J)$ is called  {\it  $k$-nondegenerate}  if   $\cF_j\! \neq\! 0$ for all $0 \leq$ $j$ $\leq k-2$  and $\cF_{k - 1} \!\!=\! 0$. In this case,  $(M, \cD, J)$ is  called {\it finitely nondegenerateÊ with    order of nondegeneracy} $k$.
\end{definition}
If $(M, \cD, J)$ is of hypersurface type (i.e.  of   codimension $1$), it is $1$-nondegenerate if and only if  it is Levi nondegenerate in the usual sense or, equivalently, if and only if  $\cD$ is a contact distribution.  \par
\smallskip
If $(M, \cD, J)$ is of dimension $5$ and of hypersurface type,  a  dimension argument shows that it is finitely nondegenerate if and only if it is either Levi nondegenerate or  $2$-nondegenerate.  \par
\medskip 
In  this paper, we focus on   {\it 5-dimensional, CR manifolds  of hypersurface type that are $2$-nondegenerate}, which  we    friendly call  {\it girdled CR manifolds}. \par
\smallskip
By definition, a girdled CR manifold $(M, \cD, J)$ is naturally endowed with the $J$-invariant, 2-dimensional, involutive  subdistribution $\cE = \Re(\cF_0) \subset \cD$, which we  call     {\it rib distribution}.  Its maximal leaves  are called  {\it ribs}:  they are  complex manifolds of  dimension $1$ and  $(M, \cD, J)$ is foliated by such 1-dimensional complex manifolds. However, by  finite nondegeneracy, there exist  no  local CR equivalences  between $(M, \cD, J)$  and  products of the form $\overline M \times \bC$ for  some 3-dimensional CR manifold $(\overline M, \overline \cD, \overline J)$ (we call them {\it CR straightenings}; see also \cite{Fr}).  The absence of CR straightenings is the reason why we chose  the word   ``girdled''    for  such CR manifolds. \par
\subsection{Levi form and  cubic form}\hfill\par
Let $(M, \cD, J)$ be a girdled CR manifold with rib distribution $\cE \subset \cD$  and   denote by 
$$\cE^{10} = \cD^{10} \cap \cE^\bC\,,\qquad \cE^{01} = \overline{\cE^{10}}$$ 
 the holomorphic and antiholomorphic subdistributions in $\cE^\bC$. For any  tangent vector $v \in T_x M$, we   use the notation  $X^{(v)}$  to indicate any  vector field, defined on a neighbourhood of $x$, satisfying the condition
$$X^{(v)}|_x = v\,.$$
Any such  $X^{(v)}$  is said to be  {\it a vector  field that extends $v$ around $x$}.\par
\smallskip 
We call  {\it  defining covector at $x \in M$} a 1-form   $\vartheta_x \in T_{x}^* M$  with the property  that $\ker \vartheta_x = \cD_{x}$.  A 1-form $\vartheta$,  defined on an open subset  $\cU \subset M$, is called  {\it  defining 1-form}  if    $\vartheta_x$ is a defining covector for any  $x \in \cU$. \par
\begin{lem} Let $\vartheta$ be a defining 1-form on a neighbourhood $\cU$ of $x \in M$. 
\begin{itemize}
\item[a)] For any $v, w \in T_x M$ and vector fields $X^{(v)}$, $X^{(w)} \in \under \cD$ that extend $v$ and $w$ around $x$, we have that 
\beq \label{Leviform1}Êd \vartheta_x(v, J w) = - \vartheta_x([X^{(v)}, JX^{(w)}])\,.\eeq
In particular,  $d \vartheta_x(v, J w)$ depends  only on $\vartheta_x$, $v$ and $w$. 
\item[b)] Given  $e \in \cE^{10}_x$, $h, h' \in \cD^{01}_x$ and    $X^{(e)} \in \under \cE^{10}$, $X^{(h)}$, $X^{(h')} \in \under \cD^{01}$ that extend $e$, $h$ and $h'$, respectively, the corresponding complex number  
$$\vartheta_x([[X^{(e)}, X^{(h)}], X^{(h')}])$$
  depends only on $\vartheta_x$,  $e$,  $h$ and $h'$. Such dependence is linear.
\end{itemize}
\end{lem}
\begin{pfns} The equality \eqref{Leviform1}   is a consequence of Koszul formula for  exterior derivatives.  The last claim of (a)  follows directly. \par
For (b), we only need to check that  the value of $\vartheta_x([[X^{(e)}, X^{(h)}], X^{(h')}])$  does not change if one replaces $X^{(e)} \in \under \cE^{10}$, $X^{(h)}$, $X^{(h')} \in \under \cD^{01}$  by other  extensions $Y^{(e)} \in \under \cE^{10}$, $Y^{(h)}$, $Y^{(h')} \in \under \cD^{01}$.   They are necessarily of the form 
$$Y^{(e)} = \l X^{(e)}\,,\quad Y^{(h)} = \mu \overline{X^{(e)}} + \nu X^{(h)} , \quad Y^{(h')} = \mu' \overline{X^{(e)}} + \nu' X^{(h')}$$
for some $\bC$-valued, smooth  functions $\l$, $\m$, $\nu$, $\mu'$, $\nu'$ with 
$$\l_x =  \nu_x = \nu'_x = 1\,,\qquad \mu_x= \mu'_x = 0\,.$$
Since  $[\under \cD^{01},\under  \cD^{01}] \subset\under \cD^{01}$ and  $[\under\cE^\bC, \under\cD^\bC] \subset \under \cD^\bC$,   
$$[Y^{(e)}, Y^{(h)}] =  \l \nu [X^{(e)}, X^{(h)}] + \l X^{(e)}(\nu) X^{(h)} \mod\under \cE^\bC$$
and  
$$[[Y^{(e)}, Y^{(h)}] , Y^{(h')}]  =  \l \nu \nu' [[X^{(e)}, X^{(h)}], X^{(h')}]  \mod  \under \cD^{\bC}\,.$$
 Since $\vartheta$ is a defining 1-form and $\l_x = \m_x = \nu_x = 1$
$$\vartheta_x([[Y^{(e)}, Y^{(h)}] , Y^{(h')}] ) =  \vartheta_x([[X^{(e)}, X^{(h)}], X^{(h')}])\,. \eqno\qed$$ 
\end{pfns}
On the base of the previous lemma, we may consider the following\par
\begin{definition} Let $\vartheta_x \in T_{x}^* M$ be a defining covector at $x$. We call {\it  Levi form} and {\it  cubic form},  associated with $\vartheta_x$,    the linear maps
\beq \label{definitionLevi} \cL^{\vartheta_x}: \cD_x \times \cD_x \longrightarrow \bR\,,\qquad \cL^{\vartheta_x}(v, w) = - \vartheta_x([X^{(v)}, J X^{(w)}])\,,\eeq
\beq \label{definitioncubic}Ê\cH^{\vartheta_x}: \cE^{10}_x \times \cD^{01}_x \times \cD^{01}_x \longrightarrow \bC\,,\qquad  \cH^{\vartheta_x}(e, h, h') = 
\vartheta_x([[X^{(e)}, X^{(h)}], X^{(h')}])\eeq
for some  extensions $X^{(v)}, X^{(w)} \in \under \cD$,  $X^{(e)} \in \under \cE^{10}$,  $X^{(h)}$, $ X^{(h')} \in \under \cD^{01}$. 
\end{definition}
\section{A maximally homogeneous model for girdled CR manifolds: \\the tube over the future light cone}
\label{homogeneous model}
\setcounter{equation}{0}
\subsection{The tube over the future light cone}\hfill\par
Consider the bilinear form $(\cdot , \cdot )$ and the pseudo-Hermitian form  $<\cdot , \cdot>$ on $\bC^5$ defined by 
$$(t, s) = t^T I_{3,2}  s\,,\qquad <t, s> = (\overline t, s) \,,\qquad  I_{3,2} =  \left(\begin{array}{c|c} I_3 & 0\\ \hline 
0 & - I_2 \end{array} \right),$$
and the corresponding semi-algebraic subset   $M_o \subset \bC P^4$  defined by 
\beq\label{tubeflc1} \left\{\begin{array}{l}  (t, t) =  (t^0)^2 + (t^1)^2 + (t^2)^2 - (t^3)^2   -  (t^4)^2= 0\,,\\
\,\\
< t, t> =  |t^0|^2 + |t^1|^2 + |t^2|^2 - |t^3|^2   - |t^4|^2 = 0\,,\\
\, \\
\Im\left(t^3 \overline{ t^4}\right) >  0\,.
\end{array}\right.\eeq
 One can directly check (see also e.g. \cite{SV}) that  $M_o$ 
  is a  $\SO_{3,2}^o$-homogeneous,  5-dimensional CR submanifold  of $\bC P^4$  ($\SO_{3,2}^o$ = identity component of $\SO_{3,2}$)   and contains  $\wt \cT = M_o \cap \{ \Im(t^3 \overline{(t^0 + t^4)}) > 0\}$  as   open   dense  subset, which  is  CR equivalent to the
so called {\it tube over the future light cone in $\bC^3$}, i.e. the real hypersurface
\beq\label{tubeflc2} \cT = \{\,(z^1, z^2, z^3) \in \bC^3\,: \,(x^1)^2 + (x^2)^2 - (x^3)^2 = 0\,,\,x^3 > 0\,\}\,. \eeq
In fact, one can directly check that the map $f: \bC^3\longrightarrow \bC P^4$ defined by 
\beq \label{expr}Ê f(z^1, z^2, z^3) = \left[ - \frac{i}{2} - \frac{i}{2} \left( (z^1)^2 + (z^2)^2 - (z^3)^2  \right): z^1 : \phantom{aaaaaaaaaaaaaaaaa}\right.$$
$$\left.\phantom{aaaaaaaaaaaaaaaaaaaa} : z^2 :  z^3 : - \frac{i}{2} +  \frac{i}{2} \left( (z^1)^2 + (z^2)^2 - (z^3)^2  \right)\right]\eeq
determines a CR equivalence  between $\cT$ and  $\wt \cT \subset M_o$. \par
\medskip
It is also  known that  $\cT$ (and  $M_o$ as well) is a girdled CR manifold. It is indeed  a homogeneous girdled CR manifold with  algebra  of germs of infinitesimal automorphisms of maximal  dimension (see (\cite{FK1}). Indeed
 the real algebraic  variety 
$$N = \{\,[t] \in \bC P^4\,: \,(t, t)  = (\overline t, t)  = 0\,\} \subset \bC P^4\,.$$
is   $\operatorname{O}_{3,2}$-invariant    and contains exactly  two open Ê$\SO_{3,2}^o$-orbits,   one of which is $M_o$.   It is known that   $\Aut(M_o) = \SO_{3,2}^o$  , i.e. its CR automorphisms coincide with the transformations determined by the   projective actions  of the elements in 
Ê$\SO_{3,2}^o$ on $M_o$ (see e.g. \cite{FK}). 
\subsection{The graded structure of the Lie algebra of $\Aut(M_o)$}
\label{section3.2}\hfill\par
Consider  a   system of projective coordinates on $\bC P^4$, in which   the bilinear form $(\cdot, \cdot)$ assumes the form
\beq \label{projcoor}Ê(t, s) = t^T\cI\,s\qquad  \text{with}\qquad \cI = \left(
\begin{array}{cc|c|cc}0&0&0&0&1\\
0&0&0&1&0\\
\hline
0&0&1&0&0\\
\hline
0&1&0&0&0\\
1&0&0&0&0
\end{array}
\right)\,.\eeq
By means of these new coordinates,  $ \so_{3,2}$  can be identified with the  Lie algebra of  real matrices such that  $A^T \cI + \cI A = 0$,  i.e.,  of the form 
$$A =  \left(
\begin{array}{c|c|c}
\begin{matrix} a_1 & a_2 \\a_3 & a_4 \end{matrix} &  \begin{matrix}  a_5 \\a_6 \end{matrix} & \begin{matrix}  a_7 & 0 \\0 & -a_7 \end{matrix}\\
\hline 
\begin{matrix}   a_{8} &  a_{9} \end{matrix} & 0  & \begin{matrix}  -a_{6} & -a_5 \end{matrix}\\
 \hline
\begin{matrix} a_{10} & 0 \\0 & -a_{10} \end{matrix} &  \begin{matrix}  - a_{9}\\-a_{8} \end{matrix} & \begin{matrix}- a_4  & - a_2 \\- a_3 & - a_1 \end{matrix}
\end{array}
\right)\,,\qquad\text{for some}\, a_i \in \bR\,.$$
In particular,  it  admits a basis $\cB$,  given by the matrices
\beq \label{basisso32}{\tiny  
e^{-2} =  \left(
\begin{array}{c|c|c}
\begin{matrix} 0 & 0 \\0 &0 \end{matrix} & \begin{matrix} 0 \\ 0  \end{matrix} &  \begin{matrix}  0 & 0 \\0 & 0 \end{matrix}\\
\hline 
\begin{matrix}  0 &   0 \end{matrix} & 0  & \begin{matrix}  0 &   0 \end{matrix}\\
 \hline
\begin{matrix}  1 & 0 \\0 & -1 \end{matrix} &  \begin{matrix} 0 \\ 0  \end{matrix} & \begin{matrix} 0  & 0 \\0 & 0 \end{matrix}
\end{array}
\right),\quad 
e^{-1}_1 =  \left(
\begin{array}{c|c|c}
\begin{matrix} 0 & 0 \\0 &0 \end{matrix} & \begin{matrix} 0 \\ 0  \end{matrix} &  \begin{matrix}  0 & 0 \\0 & 0 \end{matrix}\\
\hline 
\begin{matrix}  1 &   0 \end{matrix} & 0  & \begin{matrix}  0 &   0 \end{matrix}\\
 \hline
\begin{matrix}  0 & 0 \\0 & 0 \end{matrix} &  \begin{matrix} 0 \\ -1  \end{matrix} & \begin{matrix}0  & 0 \\0 & 0 \end{matrix}
\end{array}
\right),} $$
$$ {\tiny \, 
e^{-1}_2 =  \left(
\begin{array}{c|c|c}
\begin{matrix} 0 & 0 \\0 &0 \end{matrix} & \begin{matrix} 0 \\ 0  \end{matrix} &  \begin{matrix}  0 & 0 \\0 & 0 \end{matrix}\\
\hline 
\begin{matrix}  0 &   1 \end{matrix} & 0  & \begin{matrix}  0 &   0 \end{matrix}\\
 \hline
\begin{matrix}  0 & 0 \\0 & 0 \end{matrix} &  \begin{matrix} -1 \\ 0  \end{matrix} & \begin{matrix}0  & 0 \\0 & 0 \end{matrix}
\end{array}
\right), 
\quad
e^0_1 =  \left(
\begin{array}{c|c|c}
\begin{matrix} 1 & 0 \\0 & - 1 \end{matrix} & \begin{matrix} 0 \\ 0  \end{matrix} &  \begin{matrix}  0 & 0 \\0 & 0 \end{matrix}\\
\hline 
\begin{matrix}  0 &   0 \end{matrix} & 0  & \begin{matrix}  0 &   0 \end{matrix}\\
 \hline
\begin{matrix}  0 & 0 \\0 & 0 \end{matrix} &  \begin{matrix} 0 \\ 0  \end{matrix} & \begin{matrix} 1  & 0 \\0 & -1 \end{matrix}
\end{array}
\right),}
$$
$$ {\tiny  
e^0_2 =  \left(
\begin{array}{c|c|c}
\begin{matrix} 0 & 1 \\1 &0 \end{matrix} & \begin{matrix} 0 \\ 0  \end{matrix} &  \begin{matrix}  0 & 0 \\0 & 0 \end{matrix}\\
\hline 
\begin{matrix}  0 &   0 \end{matrix} & 0  & \begin{matrix}  0 &   0 \end{matrix}\\
 \hline
\begin{matrix}  0 & 0 \\0 & 0 \end{matrix} &  \begin{matrix} 0 \\ 0  \end{matrix} & \begin{matrix}0  & -1 \\-1 & 0 \end{matrix}
\end{array}
\right),
\quad
E^0_1 =  \left(
\begin{array}{c|c|c}
\begin{matrix} 1 & 0 \\0 &1 \end{matrix} & \begin{matrix} 0 \\ 0  \end{matrix} &  \begin{matrix}  0 & 0 \\0 & 0 \end{matrix}\\
\hline 
\begin{matrix}  0 &   0 \end{matrix} & 0  & \begin{matrix}  0 &   0 \end{matrix}\\
 \hline
\begin{matrix}  0 & 0 \\0 & 0 \end{matrix} &  \begin{matrix} 0 \\ 0  \end{matrix} & \begin{matrix}- 1  & 0 \\0 & -1 \end{matrix}
\end{array}
\right),}
$$
$$ {\tiny 
E^0_2 =  \left(
\begin{array}{c|c|c}
\begin{matrix} 0 & 1 \\-1 &0 \end{matrix} & \begin{matrix} 0 \\ 0  \end{matrix} &  \begin{matrix}  0 & 0 \\0 & 0 \end{matrix}\\
\hline 
\begin{matrix}  0 &   0 \end{matrix} & 0  & \begin{matrix}  0 &   0 \end{matrix}\\
 \hline
\begin{matrix}  0 & 0 \\0 & 0 \end{matrix} &  \begin{matrix} 0 \\ 0  \end{matrix} & \begin{matrix}0  & -1 \\1 & 0 \end{matrix}
\end{array}
\right),
\quad
E^1_1 =  \left(
\begin{array}{c|c|c}
\begin{matrix} 0 & 0 \\0 &0 \end{matrix} & \begin{matrix} 1 \\ 0  \end{matrix} &  \begin{matrix}  0 & 0 \\0 & 0 \end{matrix}\\
\hline 
\begin{matrix}  0 &   0 \end{matrix} & 0  & \begin{matrix}  0 &   -1 \end{matrix}\\
 \hline
\begin{matrix}  0 & 0 \\0 & 0 \end{matrix} &  \begin{matrix} 0 \\ 0  \end{matrix} & \begin{matrix}0  & 0 \\0 & 0 \end{matrix}
\end{array}
\right),}
$$
$$ {\tiny  
E^1_2 =  \left(
\begin{array}{c|c|c}
\begin{matrix} 0 & 0 \\0 &0 \end{matrix} & \begin{matrix} 0 \\ 1  \end{matrix} &  \begin{matrix}  0 & 0 \\0 & 0 \end{matrix}\\
\hline 
\begin{matrix}  0 &   0 \end{matrix} & 0  & \begin{matrix}  -1 &   0 \end{matrix}\\
 \hline
\begin{matrix}  0 & 0 \\0 & 0 \end{matrix} &  \begin{matrix} 0 \\ 0  \end{matrix} & \begin{matrix}0  & 0 \\0 & 0 \end{matrix}
\end{array}
\right),
\quad
E^2 =  \left(
\begin{array}{c|c|c}
\begin{matrix} 0 & 0 \\0 &0 \end{matrix} & \begin{matrix} 0 \\ 0  \end{matrix} &  \begin{matrix}  1 & 0 \\0 & -1 \end{matrix}\\
\hline 
\begin{matrix}  0 &   0 \end{matrix} & 0  & \begin{matrix}  0 &   0 \end{matrix}\\
 \hline
\begin{matrix}  0 & 0 \\0 & 0 \end{matrix} &  \begin{matrix} 0 \\ 0  \end{matrix} & \begin{matrix} 0  & 0 \\0 & 0 \end{matrix}
\end{array}
\right).}\eeq
For the discussions of the next sections, it is  quite useful
to have  all Lie brackets between  elements of  the basis $\cB$ explicitly written down. Moreover, in place  of the   elements $E^{i}_j$ and $e^{-i}_j$,  for $i = 0,1$ and $j = 1,2$,  it is   convenient to consider  the   elements in $(\so_{3,2})^\bC$ 
\beq \label{definitionholelements}
\begin{array}{llll}
E^{\ell (10)} &=  \frac{1}{2} \left(E^\ell_1 - i E^\ell_2\right)\,, &
E^{\ell(01)} &= \overline{E^{\ell(10)}}
\,,\\
\,\\
e^{-\ell(10)} &=  \frac{1}{2} \left(e^{-\ell}_1 - i e^{-\ell}_2\right)
 \,, &e^{-\ell(01)} &= \overline{e^{-\ell(10)}}\,,
 \end{array} \quad \ell = 0,1\,,\eeq
and  evaluate the Lie brackets between   such vectors and between  them and  the real  vectors  $E^2$,  $e^{-2}$ or $E^0_1$.  Here is the list  of such Lie brackets.\par
 \smallskip
  \renewcommand{\arraystretch}{2}
\renewcommand{\arraycolsep}{2pt}
\par
 \medskip
 \centerline{\rotatebox{0}{$\tiny\begin{array}{||c|c|c|c|c|c|c|c|c|c|c||}\hline\hline \,&
E^2& 
E^{1(10)}&
E^{1(01)}&
E^{0(10)}&
E^{0(01)}&
e^{0(10)}&
e^{0(01)}& 
e^{-1(10)}& 
e^{-1(01)}&
e^{-2} 
\\
\hline
\hline
\hline
\ad_{E^0_1} &
2 E^2& 
E^{1(10)}&
E^{1(01)}&
0 &
0&
0&
0& 
- e^{-1(10)}& 
- e^{-1(01)}&
- 2 e^{-2} 
\\
\hline
\hline
\hline
\ad_{E^{2}} &
{\bf 0} & 
0 & 
 0 &
-E^2 &
- E^2 &
0 &
0 &
 i E^{1(10)}  &
- i E^{1(01)} &
\stackrel{\textstyle E^{0(10)}Ê+}{+ E^{0(01)}}
\\
\hline
\ad_{E^{1(10)}} &
0 & 
 {\bf 0} & 
- \frac{i}{2} E^2 &
0 &
- E^{1(10)} &
0 &
- E^{1(01)} &
\frac{1}{2} e^{0(10)} &
\frac{1}{2} E^{0(01)} &
i e^{-1(10)}
\\
\hline
\ad_{E^{1(01)}} &
0 & 
 \frac{i}{2} E^2  & 
{\bf 0} &
- E^{1(01)} &
0  &
- E^{1(10)} &
0 &
\frac{1}{2} E^{0(10)}&
\frac{1}{2} e^{0(01)} &
- i e^{-1(01)}
\\
\hline
\ad_{E^{0(10)}} &
E^2  & 
0  & 
E^{1(01)} &
{\bf 0} &
0 &
- e^{0(10)} &
e^{0(01)} &
- e^{-1(10)} &
0 &
- e^{-2}
\\
\hline
\ad_{E^{0(01)}} &
E^2 & 
 E^{1(10)} & 
 0 &
0 &
{\bf 0} &
e^{0(10)} &
- e^{0(01)} &
0 &
- e^{-1(01)} &
-e^{-2} \\
\hline
\ad_{e^{0(10)}} &
0  & 
0 & 
E^{1(10)}&
e^{0(10)} &
-e^{0(10)} &
{\bf 0} &
\stackrel{\textstyle{ - E^{0(10)}+}}{+ E^{0(01)}} &
0 &
- e^{-1(10)} &
 0\\
\hline
\ad_{e^{0(01)}} &
0   & 
 E^{1(01)} & 
0 &
-e^{0(01)} &
e^{0(01)}  &
\stackrel{\textstyle{ E^{0(10)}Ê-}}{- E^{0(01)}}  &
{\bf 0} &
- e^{-1(01)} &
0 &
 0 \\
\hline
\ad_{e^{-1(10)}} &
- i E^{1(10)}  & 
 -\frac{1}{2} e^{0(10)} & 
 -\frac{1}{2} E^{0(10)}  &
e^{-1(10)} &
0 &
0 &
e^{-1(01)}  &
{\bf 0} &
\frac{i}{2} e^{-2}  &
0 \\
\hline
\ad_{e^{-1(01)}} &
i E^{1(01)} & 
-  \frac{1}{2} E^{0(01)}& 
-  \frac{1}{2} e^{0(01)} &
0 &
e^{-1(01)} &
e^{-1(10)} &
0 &
- \frac{i}{2} e^{-2} &
{\bf 0} &
 0 \\
\hline
\ad_{e^{-2}} &
\stackrel{ \textstyle{-  E^{0(10)}-}}{Ê- E^{0(01)}}& 
-  i e^{-1(10)} & 
i e^{-1(01)} &
e^{-2} &
e^{-2} &
0 &
0 &
0 &
0  &
{\bf 0}
\\
\hline
\hline
\end{array}$
}}\par
\medskip
\centerline{\small \bf Table 1}
\par\noindent
The  brackets between all the elements of  $\cB$ can be  directly recovered from Table 1  recalling that
\beq  \begin{array}{llll}
E^\ell_1 &= E^{\ell(10)} + E^{\ell(01)}\,, &E^\ell_2 &= i (E^{\ell(10)} - E^{\ell(01)})\,,\\
e^{-\ell}_1 &= e^{-\ell(10)} + e^{-\ell(01)} \,, &e^{-\ell}_2 &= i (e^{-\ell(10)} - e^{-\ell(01)})\,,
 \end{array} \qquad \,\ell = 0,1\,.\eeq
Notice that  $\so_{3,2}$ has a natural  graded Lie algebra structure 
$$\so_{3,2} = \gg^{-2} + \gg^{-1} + \gg^0 + \gg^1 + \gg^2\,,$$
given by the eigenspaces of the adjoint action of the {\it grading element} $E^0_1$
$$\gg^{-2} = < e^{-2} >\,,\,\gg^{-1} = < e^{-1}_1, e^{-1}_2>\,,\,\gg^0 = < e^0_1, e^0_2, E^0_1, E^0_2> \,,$$
$$\gg^1 = < E^1_1, E^1_2>\,,\,\gg^2 = < E^2>\,.$$
\subsection{The Levi form and  the cubic form of $M_o$}\hfill\par
\label{Levicubicsection}
Consider the point $x_o = [1 : i : 0: 0: 0]$ of $M_o$ (we are  using the  coordinates defined   by \eqref{projcoor}). The isotropy subalgebra $\gh = Lie (H)$ with $H \=\Aut_{x_o}(M_o)$, consists of  the matrices $A \in \so_{3,2}$ such that
$$A \cdot (1,  i ,  0,  0, 0) = \l (1,  i ,  0,  0, 0)\qquad \text{for some} \,0 \neq \l \in \bC\,.$$ 
From this, one can directly check that 
 $$\gh = < E^0_1, E^0_2, E^1_1, E^1_2, E^2 > = \gh^0 + \gg^1 + \gg^2\,,$$
 where $\gh^0 \subset \gg^0$ is the subspace   $\gh^0 = < E^0_1, E^0_2>$.
If we denote by  $\gm^0 \subset \gg^0$ the subspace
$\gm^0 = < e^0_1, e^0_2> $, 
we see that   $\gm = \gg^{-2} + \gg^{-1}Ê+ \gm^0$ is a complementary subspace  to $\gh$ in $\gg$ and  that the linear map 
$$\imath: \gm\longrightarrow T_{x_o} M\,,\qquad X \in \gm\quad  \overset\imath\longmapsto\quad  \under X \= \wh X|_{x_o} \in T_{x_o} M_o$$
 where  $\wh X$ is the  vector field on $M_o$  with flow equal to $\exp(t X)\cdot: M_o \longrightarrow M_o$, $t \in \bR$, is  an isomorphism. We    use $\imath$   to identify   $\gm$ with $T_{x_o} M$ and,  for any $v \in \gm$, we  set $\underline v = \imath(v)$. \par
In this way,  the basis $\cB^\gm = (e^{-2}, e^{-1}_i, e^0_i)$ of $\gm$ can be identified with the  basis  $\under \cB^\gm$ of   $T_{x_o} M_o$ formed  by the vectors 
$$\underline e^{-2} = \Re\left(\left.\frac{\partial}{\partial z^3}\right|_{x_o} - i \left.\frac{\partial}{\partial z^4}\right|_{x_o}\right) \,,\,$$
$$ \under e^{-1}_1  = \Re\left(\left.\frac{\partial}{\partial z^2}\right|_{x_o}\right)\,,\qquad  \under e^{-1}_2 = \Re\left(\left.i \frac{\partial}{\partial z^2}\right|_{x_o}\right)\,,$$
$$\under e^{0}_1  = \Re\left(\left.\frac{\partial}{\partial z^0}\right|_{x_o} - i\left. \frac{\partial}{\partial z^1}\right|_{x_o}\right)\,,\qquad \under  e^{0}_2  =\Re\left(\left. i \frac{\partial}{\partial z^0} \right|_{x_o}+  \left.\frac{\partial}{\partial z^1}\right|_{x_o}\right)\,. $$
From this  and with the help  of Table 1, one can   check that  the invariant CR  distribution $\cD$ and the rib distribution $\cE$ are  given  by
\beq \label{3.8} \cD_{x_o} = <\under e^{-1}_1,  \under e^{-1}_2, \under e^{0}_1, \under e^{0}_2 >  = \imath(\gg^{-1}Ê+ \gm^0)\,,\,\, \cE_{x_o} = < \under e^{0}_1, \under e^{0}_2 > = \imath(\gm^0)\,,\eeq
and    the invariant complex structure  $J$ on $\cD$  is such that 
$$J (\under e^{-1}_1) = \under e^{-1}_2\,,\qquad J(\under e^0_1) = \under e^0_2\,.$$
Moreover, if we denote by $\under \cB^{\gm*}$ the dual basis of $\under \cB^\gm$, direct computations show that 
$$\vartheta = (\under e^{-2})^*\qquad (\,= \text{the 1-form dual to} \, \under e^{-2}\,)$$
  is a 
defining covector for $\cD$, with  corresponding  Levi form and cubic form   equal to 
\beq \label{LeviformofTLC}Ê\cL^\vartheta(\under e^{-1(10)}, \under e^{-1(01)}) = - \frac{1}{2}\,,\qquad \cH^\vartheta(\under e^{0(10)}, \under e^{-1(01)}, \under e^{-1(01)}) = - \frac{i}{2} \eeq
(here,   $\under e^{-\ell (10)} = \frac{1}{2}\left(\under e^{-\ell}_1 - i J \under e^{-\ell}_1\right)$ and $\under e^{-\ell (01)} = \overline{\under e^{-\ell(10)}}$). 
In fact,   consider the  structure of $H$-principle bundle of 
  $\SO_{3,2}^o = \Aut(M_o)$ onto $M_o$
$$\pi:  \SO_{3,2}^o  \longrightarrow  \SO_{3,2}^o/ H \simeq M_o\ ,\qquad 
H = \Aut_{x_o}(M_o)$$
and observe  that,  by  \eqref{3.8},  the CR-structure of $ M_o$ is 
formed by  the  distribution, spanned  by the vectors $\pi_*(e^{-1}_i 
|_a)$,
  $\pi_*(e^0_j |_a)$, $i, j = 1,2$, $a \in  \SO_{3,2}^o$, and  the 
family of complex structures
  $J(\pi_*(e^i_1 |_a)) = \pi_*(e^i_2 |_a)$,   $i = -1$, $0$.  Hence, for 
a fixed local section $\s: \cU \subset M_o \longrightarrow \SO_{3,2}^o$ 
about the origin,  the vector fields $X \in \under \cD^{10}|_{\cU}$  are 
$\s$-related to vector fields of the form $\wh X =  \l_{-1} e^{-1(10)} + 
\l_{0} e^{0(10)}\!\!\! \mod \gh$ on $\SO_{3,2}^o$,  for some  
$\bC$-valued functions  $\l_{-1}$, $\l_0$. Using this and definitions, 
one can  check  \eqref{LeviformofTLC} using Table 1.  \par
\subsection{The CR structure of $M_o$ and the    complex structure $J|_{\gm^{-1} + \gm^0}$}\hfill\par
\label{complexstructureofgg}
The graded subspaces $\gg^{-2}$, $\gg^{-1}$, $\gg^1$, $\gg^2\subset \so_{3,2}$ will be often  indicated  with the symbols 
$\gm^{-2}$, $\gm^{-1}$, $\gh^1$, $\gh^2$, respectively,  so that   the graded decomposition of $\gm$ and $\gh$ are   
$$\gm = \gm^{-2} + \gm^{-1} + \gm^0\qquad \text{and}\qquad  
\gh = \gh^0 + \gh^1 + \gh^2\,.$$ 
We  also denote by $J$ the grade preserving  complex structure on the subspace 
$\gm^{-1} + \gm^0 + \gh^0 + \gh^1$,  defined by 
\beq \label{cstrgg} J(e^{-1}_1) =  e^{-1}_2 \,,\quad J(e^{0}_1) =  e^{0}_2\,,\quad J(E^0_1) = E^0_2\,,\quad J(E^1_1) = E^1_2\,.\eeq
Notice that: 
\begin{itemize}
\item[--] the  vectors $e^{-\ell(10)}$, $E^{\ell(10)}$ and $e^{-\ell(01)}$, $E^{\ell(01)}$,   introduced in \eqref{definitionholelements}, are  the $J$-holomorphic and $J$-antiholomorphic parts of the elements $e^{-\ell}_1$, $E^{\ell}_1$; 
\item[--]  through the  isomorphism $\imath: \gm \longrightarrow  T_{x_o}ÊM_o$, {\it the restriction $J|_{\gm^{-1} + \gm^0}$ corresponds to   the complex structure $J$ of $\cD_{x_o}$} and completely determines the invariant CR structure of $M_o$; 
\item[--] using  the isomorphism $\imath$, the  vectors $e^{-\ell(10)}$, $\ell = 0,1$,  can be considered as  a basis for  $\cD^{10}_{x_o} \subset T^\bC_{x_o} M_o$.
\end{itemize}
\section{Filtered vector spaces  modelled on  $\so_{3,2}$}
\setcounter{equation}{0}
\subsection{Filtrations  with an additional semitone}\hfill\par
In the following sections,    we  have   to study the   geometric structures on manifolds modelled on  
the Lie algebra $\gg = \so_{3,2}$.  In order to do this,   we  want   to establish a  few   properties of 
filtrations of vector spaces that are similar to some special  filtrations of    $\so_{3,2}$. 
\par
\medskip
Consider a finite-dimensional vector space  $V$ endowed with a filtration   of the following kind
\beq \label{filtration} \cF\,:\,V = V_{-2} \supset V_{-1} \supset V_{0} \supset V_1 \supset \ldots \supset V_k \supset \{0\}\,.\eeq
If   $V_{0}$ has its own filtration with an extra term  $V_{(0|0)}$ between $V_0$ and  $V_1$, i.e. 
\beq \label{filtrationV0}  \underset{=V_{0}}{V_{(0|-1)}}Ê\supsetneq V_{(0|0)} \supset  \underset{=V_{1}}{ V_{(0|1)}}\supset \ldots \supset  \underset{=V_{k}}{V_{(0|k)}} \supset  \{0\}\,,\eeq
 we may  merge such two filtrations and obtain a new one, namely
\beq \label{filtrationstar} \cF_* \,:\,V = V_{-2} \supset V_{-1} \supset V_{0}\supsetneq	 V_{(0|0)}  \supset V_1 \supset \ldots \supset  V_k \supset   \{0\}\,.\eeq
Such new filtration  is  called   {\it filtration with one additional  semitone}. The main example to have in  mind is given by the two filtrations
of $\gg = \so_{3,2}$ 
\beq \label{filtrationgg} \cF \,: \,\underset{V_{-2}}{\underbrace{\gm + \gh}} \supset \underset{V_{-1}}{\underbrace{\gm^{-1} + (\gm^0 + \gh^0 + \gh^1 + \gh^2)}}\supset \underset{V_0}{\underbrace{ \gm^0+ \gh^0 + \gh^1 + \gh^2}} \supset\phantom{aaaaaaa}$$
$$\phantom{aaaaaaaaaaaaaaaaaaaaaaaaaaaaaaaaa} \supset
\underset{ V_1}{\underbrace{ \gh^1 + \gh^2} }\supset \underset{V_2}{\underbrace{ \gh^2}} \supset \{0\}\,,\eeq
\beq \label{filtrationstargg}Ê\cF_* \,: \\\underset{V_{-2}}{\underbrace{\gm +\gh}} \supset \underset{V_{-1}}{\underbrace{\gm^{-1} + (\gm^0 + \gh^0 + \gh^1 + \gh^2)}}\supset \underset{V_0 = V_{(0|-1)}}{\underbrace{ \gm^0+ \gh^0 + \gh^1 + \gh^2}} \supset\phantom{aaaaaaa}$$
$$\phantom{aaaaaaaaaaaaaaaaaaa}\supset \underset{V_{(0|0)}}{\underbrace{\gh^0 + \gh^1 + \gh^2}} \supset
\underset{ V_{(0|1)}Ê= V_1}{\underbrace{ \gh^1 + \gh^2} }\supset \underset{V_{(0|2)} = V_2}{\underbrace{ \gh^2}} \supset \{0\}\,.\eeq
In the following, given  filtrations \eqref{filtration} and \eqref{filtrationstar}, we  denote by  $\GL(V, \cF)$, $\GL(V, \cF_*)$ and   $\ggl(V, \cF)$, $\ggl(V, \cF_*)$ the  Lie groups and  Lie algebras of   filtration preserving automorphisms of $(V, \cF)$ and $(V, \cF_*)$, respectively. Notice that 
 $\ggl(V, \cF)$ and $\ggl(V, \cF_*)$ are naturally endowed with  structures of  filtered Lie algebras,  with  filtrations determined by the subspaces 
$$\ggl_i(V) = \{\,A \in \ggl(V, \cF)\,:\,A(V_k) \subset V_{k+i} \,\text{ for any }\,k \geq -2\,\}\,,\phantom{aaaaaaaaa}$$
$$\ggl_{i*}(V) = \{\,A \in \ggl(V, \cF_*)\,:\, A(V_{-2})\subset V_{-2 + i}\,,\,A(V_{-1})\subset V_{-1 + i}\,,\phantom{aaaaa}$$
$$\phantom{aaaaaaaaaaaaaaaaaaaaaa} A(V_{(0|k)}) \subset V_{(0|k+i)} \,\text{for }\,k \geq -1\,\}\,, \qquad i \geq 0\,.$$
We denote by  $\GL_i(V)$ and $\GL_{i*}(V)$ the corresponding connected subgroups of $\GL(V, \cF)$ and $\GL(V, \cF_*)$, respectively.\par
\bigskip
Let  $W = \ggr(V, \cF)$ and $W_* = \ggr(V, \cF_*)$ be the  graded vector spaces of $(V, \cF)$ and $(V, \cF_*)$,  i.e. 
\beq \label{W} W =   W^{-2} + W^{-1} + W^0 + W^1 + \dots + W^k\,,\eeq
\beq \label{W_*} W_* =  W^{-2}  + W^{-1} + 
 W^{(0|-1)} + W^{(0|0)} + W^1 + \dots + W^k\eeq
with  $W^i = V_i/V_{i + 1}$ and $W^{(0|j)} = V_{(0|j)}/V_{(0|j+1)}$. The corresponding natural projections are denoted   by
  $$\pi^i: V_i \longrightarrow  W^i = V_i/V_{i+1}\,,\qquad \pi^{(0|j)}: V_{(0|j)} \longrightarrow W^{(0|j)}Ê= V_{(0|j)}/V_{(0|j +1)}\,.$$
Note  that the
graded vector spaces $W$, $W_*$  are naturally endowed with   filtrations, which we denote by $\cF$ and $\cF_*$, respectively  (the filtration $\cF$ is     $\cF = \{W_j = \sum_{i\geq j} W^i\}$; the filtration $\cF_*$ is defined analogously), so that also  the Lie groups $\GL_i(W) $ and $\GL_{i*}(W_*) $ are well defined.\par
\subsection{Partial complex structures on  filtered vector spaces}\hfill\par
From  \S \ref{complexstructureofgg}, we know  that there exist  two complex structures on the subspaces $\gm^{-1}$, $\gm^0$ of $\gg = \so_{3,2}$, which are 
algebraic counterparts of  the CR structure of $M_o$.  Motivated by this, we  consider the following
\par
\begin{definition} \label{partialJ}
  A {\it  partial complex structure on $(V, \cF_*)$} is a 
 filtration preserving linear map   $J: V_{-1}\longrightarrow V_{-1}$ such that
 \beq \label{partialJ1}ÊJ^2 = - \Id_{V_{-1}}\!\!\!\! \mod V_{(0|0)}\,.\eeq
Two partial complex structures $J$, $J'$ are called {\it equivalent} if 
$$J - J'  = 0\!\!\!\!  \mod  V_{(0|0)} \,.$$
\end{definition}
A partial complex structure $J$   induces     on $W^{-1} + W^{(0|-1)}$ the complex structure defined by, for any $X \in V_{-1}$
and $Y \in V_{(0|-1)}$,
 $$J(X_{\!\!\!\!\!\!\mod V_0}) = J (X)_{\!\!\!\!\!\!\mod V_0}\,,\qquad J(Y_{\!\!\!\!\!\!\mod V_{(0|0)}} )= J (Y)_{\!\!\!\!\!\!\mod V_{(0|0)}}\,.$$
Note that:
 \begin{itemize}
 \item[--]  equivalent partial structures   induce the same complex structure on $W^{-1} + W^{(0|-1)}$; 
 \item[--]  if  $J$ is extended to an endomorphism of  $W_{-1}$, it  
  is a partial complex structure on  $ W_* = \ggr(V, \cF_*)$; any two such extensions are equivalent. 
 \end{itemize}
 \smallskip
We  use the notation:  
$$
\begin{array}{lcl}
\GL(V, \cF_*, J) &=& \left\{\,A \in \GL(V, \cF_*)\,:  ÊJ \circ A|_{V_{-1}} = A \circ J\!\!\!\! \mod V_{(0|0)}\,\right\}\,,\\
\GL_{i*} (V , J) &=& \GL(V, \cF_*, J) \cap \GL_{i*}(V) \,,\\
\GL_{i} (V , J) &=& \GL(V, \cF_*, J) \cap \GL_{i}(V)\,,\\
\GL^{\text{gr}}_i(W, J) & =& \GL_i(W, J) \!\!\! \mod \GL_{i+1}(W, J)\, .
\end{array}
$$
Lie algebras of such groups are denoted by the corresponding gothic letters. 
\par
Note that, when $i \geq 1$, the group $\GL^{\text{gr}}_i(W, J)$ is abelian and its  Lie algebra is identifiable with the vector space  
$$
\begin{array}{lcl}
\ggl^{\text{gr}}_i(W, J) &=& \left\{\phantom{ W_{(0|0)}}\!\!\!\!\!\!\!\! \!\!\!\!\!B \in \ggl(W)\,:  Ê\, B(W^\ell)\subset W^{\ell+i} \,
 \,\text{for any}\,\ell  \,\text{and}  \right. \\
&&\phantom{aaaaaaaaaaaaa}  \left.J\circ B|_{W^{-1} + W^{(0|-1)}} = B \circ J \!\! \mod W_{(0|0)}\, \right\}\,.
\end{array}$$
\begin{rem} \label{triviality} By definition, for any $i \geq 2$, the groups $\GL_i(V, J)$ and $\GL_{i*}(V, J)$ coincide
with  the groups  $\GL_i(V)$ and  $\GL_{i*}(V)$, respectively. 
\end{rem}
 \smallskip
 Given two vector spaces $V$, $V'$, both with filtrations 
 \eqref{filtrationstar}  and  partial complex structures $J$, $J'$, respectively, 
 we call {\it filtered $(J, J')$-isomorphism} any filtration preserving isomorphism $u: V \longrightarrow V'$ such that 
 $$u \circ J = J' \circ u|_{V_{-1}} \!\!\!\! \mod V'_{(0|0)}\,.$$
Similarly, given two graded vector spaces $W_*$, $W_*'$,  with gradations  \eqref{W_*} and partial complex structures $J$, $J'$, we  call 
{\it graded  $(J, J')$-isomorphism} a grading preserving  isomorphism $\varphi: W \longrightarrow W'$ such that 
$$\varphi \circ J|_{W^{-1}Ê+ W^{(0|-1)}}  = J' \circ \varphi|_{W^{-1}Ê+ W^{(0|-1)}}\,.$$
\par
\smallskip
 Consider now the  Lie algebra $\gg = \so_{3,2}$, endowed with the filtrations \eqref{filtrationgg} and \eqref{filtrationstargg} and the complex structure 
 \eqref{cstrgg} on $\gm^{-1} +  \gm^0 + \gh^0 + \gh^1$.  {\it Any linear extension of $J$  onto  the entire subspace $\gm^{-1} +  \gm^0 +\gh$ is a partial complex structure of $\gg$} and two such extensions are equivalent.  This property  is indeed the main motivation  for the definitions considered in this section.  \par
 \smallskip
We  conclude observing that,    on the subspace $\gm = \gm^{-2}Ê+ \gm^{-1}Ê+ \gm^0 \subset \so_{3,2}$, one has that $\cF^* = \cF$ so that  $\GL_{1*}(\gm, J) = \GL_{1}(\gm, J)$.\par
 \subsection{Adapted  frames  on  spaces with partial complex structure}\hfill\par
 \label{section4.3}
In the next lemma, $V$ is a vector space with  filtrations \eqref{filtration}, \eqref{filtrationstar}  and  with a partial complex structure $J$. The  associated graded vector spaces  are $W = \ggr(V, \cF)$, $W_* = \ggr(V, \cF_*)$ and  we denote by $\wt \gm$  a (modelling) graded vector space 
$$\wt \gm = \underset{= \wt \gm^{(0|-3)}}{\wt \gm^{-2}} +   \underset{= \wt \gm^{(0|-2)}}{\wt \gm^{-1}}+  \wt \gm^{(0|-1)} +  \wt \gm^{(0|0)} +   \underset{= \wt \gm^{(0|1)}}{\wt \gm^1} + \ldots +  \underset{= \wt \gm^{(0|k)}}{\wt \gm^k}\,,$$
which is isomorphic to $W_*$ as graded vector space.  The subspace $ \wt \gm^{-1} +  \wt \gm^{(0|-1)}$ is assumed to be endowed  with a graded complex structure $\wt J$, which we 
consider extended to some graded endomorphism  of  $\wt \gm_{-1} = \sum_{j \geq -1} \wt \gm^j$, so that it can be considered as 
a partial complex structure on   $(\wt \gm, \cF_*)$, with filtration $\cF_* = \{\, \wt \gm_{(0|j)}  = \sum_{s \geq j}Ê\wt \gm^{(0|s)}\,\}$.
\par
\begin{lem} \label{lem4.2}  For any graded $(\wt J, J)$-isomorphism $u: \wt \gm \longrightarrow W  $, there exists a  filtered $(\wt J, J)$-isomorphism
$u_\sharp: \wt \gm  \longrightarrow  V$ satisfying the condition
\beq \label{usharp}Êu_\sharp(\wt \gm_{(0|j)}) \subset V_{(0|j)}\qquad \text{and}\qquad  \,u |_{\wt \gm^{(0|j)}}  = \pi^{(0|j)}\circ u_\sharp|_{\wt \gm^{(0|j)}}\eeq
for any $j \geq -3$. 
Any two such isomorphisms $u_\sharp$, $u_\sharp': \wt  \gm \longrightarrow V$ are related
by an element of $\GL_{1*}(\wt  \gm, \wt J)$, i.e., 
$$u_\sharp{}' =  u_\sharp \circ A\qquad \text{for some}\,A \in \GL_{1*}(\wt  \gm, \wt J)\,.$$
\end{lem}
\begin{pf}   For each $i < k$, consider a subspace $K^{(0|i)}$   of $V_{(0|i)}$ that is complementary to $V_{(0|i+1)}$. Assume also that the spaces $K^ {(0|-2)}$ and $K^ {(0|-1)}$ are $J$-invariant, modulo elements in $V_{(0|0)}$. The  map
$\pi = \sum_{j \geq -3} \pi^{(0|j)}$ determines a graded isomorphism between
the graded vector spaces 
$$V = K^ {(0|-3)} + \dots + K^ {(0|k-1)} + V_{(0|k)}$$
and  $W = \ggr(V)$, and   $u_\sharp  = \pi^{-1} \circ u$ is the desired isomorphism.
The last claim is immediate.
\end{pf}
The proof of Lemma \ref{lem4.2}  shows that, for a fixed graded $(\wt J, J)$-isomorphism $u: \wt \gm \longrightarrow W  $
there is a natural one to one correspondence between
 isomorphisms $u_\sharp: \wt \gm \longrightarrow V$  satisfying \eqref{usharp} and  
 ordered  sequences 
 \beq \label{sequenceH} K = (\underset{= K^ {(0|-3)}}{K^ {-2}}, \underset{= K^ {(0|-2)}}{K^ {-1}}, K^ {(0|-1)}, K^ {(0|0)},\underset{= K^ {(0|1)}}{K^ {1}}, \ldots, \underset{= K^ {(0|k-1)}}{K^ {k-1}})\,,\eeq
formed by subspaces $K^ {(0|i)} \subset V_{(0|i)}$ that are complementary to the spaces  $V_{(0|i+1)}$ and such that $K^ {(0|-2)}$ and $K^ {(0|-1)}$ are $J$-invariant modulo $V_{(0|0)}$. \par
\medskip
Any such sequence  \eqref{sequenceH} is called {\it sequence of $\cF_*$-horizontal subspaces of $V$} or, shortly, {\it s.h.s. of $V$}.  For  a fixed graded $(\wt J, J)$-isomorphism 
$u: \wt \gm \rightarrow  W$, the corresponding 
$(\wt J, J)$-isomorphism,  determined by a given choice of s.h.s. $K$,  is called {\it linear frame  of $V$  associated with  $u$ and $K$} and denoted by 
$$u_K: \wt \gm \longrightarrow V$$
or  simply $\wh K: \wt \gm \longrightarrow V$, in case  $u$ is considered as known and fixed. The  linear frames constructed in this way are called {\it adapted to the filtration and the partial complex structure  of $V$},   or  {\it  adapted}  for short.
\par
\section{The Tanaka structure of a girdled CR manifold}
\label{section2}
\setcounter{equation}{0}
\subsection{Adapted frames of a girdled CR manifold}\hfill\par
From now on,  $(M, \cD, J)$  is  a  girdled CR manifold with  rib distribution    $\cE$. 
Any tangent space $T_x M$ of $M$ is naturally endowed with a filtration $\cF$ of the form \eqref{filtration}, namely
\beq  \label{5.1}
\underset{ = V_{-2}}{T_x M} \supset \underset{ = V_{-1}}{\cD_x }\supset \underset{= V_0}{\cE_x}\supset \{0\}\,,\eeq
with  graded vector space $\ggr(T_x M, \cF)$ isomorphic to the graded subspace $\gm = \gm^{-2} + \gm^{-1} + \gm^{0}$ of $\gg = \so_{3,2}$.   Considering the  filtration $\cF_*$ of $T_x M$ of type \eqref{filtrationstar} 
 with $V_{(0|0)} = \{0\}$ 
 (and hence  with  $\cF =\cF_*$) 
 and  the complex structures $J_x: \cD_x \longrightarrow \cD_x$ of  the CR structure, we see that the $T_x M$'s  are naturally  endowed with partial complex structures. We may therefore consider the class $\cF r_0(M)$ of  $(J, J_x)$-isomorphisms 
$$u_\sharp: \gm = \gm^{-2}Ê+ \gm^{-1}Ê+ \gm^0 \longrightarrow T_x M\,,\qquad x \in M\,,$$
adapted to such filtrations and partial complex structures, i.e., 
 given  by  some graded $(J, J_x)$-isomorphism  $u: \gm \longrightarrow \ggr(T_x M)$ and  a s.h.s. $K  \subset T_x M$. \par
 \smallskip
 If  $\pi_0: \cF r_0(M) \longrightarrow M$ is the natural projection,  from  Lemma \ref{lem4.2} it follows   that  $\pi_0: \cF r_0(M) \longrightarrow M$ is a principal  bundle  with structure group 
 $$G^0_\sharp = \GL^{\text{gr}}_0(\gm, J) \ltimes \GL_{1}(\gm, J) = \GL^{\text{gr}}_0(\gm, J) \ltimes \GL_{1*}(\gm, J)\,.$$
 Note that any adapted linear frame $u_\sharp: \gm \longrightarrow T_x M$ is  uniquely determined  by the corresponding  frame  $(f^{-2}, f^{-1}_1, f^{-1}_2, f^0_1, f^0_2)$ of  $T_x M$ with
$$f^{-2} = u_\sharp(e^{-2}) \,,\qquad f^{-1}_j = u_\sharp(e^{-1}_j)\,,\qquad f^{0}_j = u_\sharp(e^{0}_j)\,,$$
or  by its  dual  coframe $(f^{-2*}, f^{-1*}_1, f^{-1*}_2, f^{0*}_1, f^{0*}_2) \subset T^*_x M$, for which 
$$\cD_x = \ker f^{-2*}\,, \qquad \cE_x = \ker f^{-2*} \cap \ker f^{-1*}_1 \cap \ker f^{-1*}_2\,.$$
In the following, for any $u_\sharp \in \cF r_0(M)$, we denote by  $\cL^u$ and $\cH^u$  the  Levi form and the   cubic form determined by the defining covector $ f^{-2*}$.
\par
\subsection{Strongly adapted frames of $(M, \cD, J)$}\hfill\par
\begin{definition}  A {\it strongly adapted frame of $T_x M$} is an adapted frame $u_\sharp = (f^{-2}, f^{-1}_1,  \ldots): \gm \longrightarrow T_x M$
such that  (compare with \eqref{LeviformofTLC})
\beq \label{adapted0}Ê\cL^u(f^{-1(10)},\overline{f^{-1(10)}}) = - \frac{1}{2}\,,\qquad \cH^u(f^{0(10)},  \overline{ f^{-1(10)}}, \overline{f^{-1(10)}}) = - \frac{i}{2}\,,\eeq
 where we set  $f^{0(10)}Ê= u(e^{0(10)})$,  $f^{-1(10)}Ê= u(e^{-1(10)})$.
\end{definition}
The following lemma is a direct consequence of definitions.
\begin{lem} The subset  $P^0_\sharp \subset \cF r_0(M)$ of  strongly adapted frames is a reduction of $\cF r_0(M)$ with  structure 
group 
$$\wt H^0_\sharp  = \wt H^0 \ltimes \GL_1(\gm,  J)\,,$$
where $\wt H^0 \subset \GL^{\mathrm{gr}}_0(\gm, J)$ is   the subgroup of maps $A$  such that
\beq\label{str1} [A(X), \overline{A(Y)}]Ê= A\left( [X, \overline Y]Ê\right)\,,\,\,[[A(Z), \overline{A(X)}], \overline{A(Y)}]Ê= A\left([[Z, \overline{X}], \overline{Y}] \right)\eeq
for any $X, Y \in \gm^{-1(10)}$,  $Z \in \gm^{0(10)}$.
\end{lem}
Since $\GL_1(\gm, J) $ is normal in $\wt H^0_\sharp$, we may consider  the quotient bundle
$$\pi^0 : P^0 = P^0_\sharp/\GL_1(\gm, J) \longrightarrow M\,,$$
which  is  a principal $\wt H^0$-bundle. Motivated by Tanaka's theory ( \cite{Ta1, Ta2, Ta3, AS}),  we call   it    {\it Tanaka structure of $(M, \cD, J)$}.\par
 \begin{rem} 
 \label{P^0interpretation}ÊWe recall that any adapted  frame $u_\sharp: \gm \longrightarrow T_x M$ is uniquely  determined by  the induced    isomorphism  $u: \gm \longrightarrow \ggr(T_x M, \cF)$ and  a  s.h.s. $K = (K^{-2}, K^{-1})$ of  $T_x M$.  Note that two adapted frames $u_\sharp$, $u'_\sharp$ are in the same   equivalence class $[u_\sharp] \in P^0|_x$ if and only if they  determine  the same graded isomorphism $u$. Hence,  $P^0$ can be also defined as  the {\it bundle of $J$-preserving, graded  isomorphisms $u: \gm \longrightarrow  \ggr(T_x M, \cF)$, determined   by    frames satisfying  
\eqref{adapted0}}. 
 \end{rem}
\begin{lem} \label{lemma5.4} The Lie algebra of $\wt H^0$  is 2-dimensional and equal to 
$$Lie(\wt H^0) = \ad(\gh^0)|_\gm\,,$$
where   $\ad(\gh^0)|_\gm\subset \ggl_0(\gm, J)$  denotes the  subalgebra  of the restrictions 
$\ad(X)|_{\gm}$,   $X \in  \gh^0 \subset \so_{3,2}$. In particular, the Lie algebra $Lie(\wt H^0)$ is abelian and  isomorphic to
$\gh^0$, and  it is naturally  endowed with the complex structure  $J|_{\gh^0}$, defined  in \eqref{cstrgg}.    
 Moreover,  $\wt H^0$ is isomorphic with the connected subgroup $H^0$ of $H \subset \SO_{3,2}^o$ with $Lie(H^0) = \gh^0$. 
\end{lem}
\begin{pf}  A   map $B \in \ggl^{\text{gr}}_0(\gm, J)$ is of the form
$$B(e^{-2}) = \t e^{-2}\,,\,B(e^{-1(10)}) = \l e^{-1(10)} = \overline{B(e^{-1(01)})}\,,\,B(e^{0(10)}) = \m e^{0(10)}$$
for some $\t \in \bR$, $\l, \m \in \bC$. On the other hand,  $B$ is in $Lie(H^0)$ if and only if 
\beq\label{str3bis} [B(e^{-1(10)}), e^{-1(01)}] + Ê [e^{-1(10)},B(e^{-1(01)})]  = B(  [e^{-1(10)}, e^{-1(01)}]) =\frac{i}{2} B(e^{-2}) \,,\eeq
\beq\label{str4bis}  [[B(e^{0(10)}), e^{-1(01)}],e^{-1(01)}]Ê+   [[e^{0(10)}, B(e^{-1(01)})], e^{-1(01)}]Ê +   $$
$$ + [[e^{0(10)}, e^{-1(01)}], B(e^{-1(01)})]Ê =   B( [[e^{0(10)} e^{-1(01)}], e^{-1(01)}]) =  \frac{-i}{2}  B(e^{-2})\,.\eeq
These conditions  are equivalent to  
$ \t = 2 \Re(\l)$ and  $\m = \t - 2 \overline \l = 2 i \Im \l$.  
Hence,  any such $B$ is  determined by the parameter $\l$ and   $Lie(H^0)$ is  spanned by 
$$B_1(e^{-2}) = - 2 e^{-2}\,,\quad  B_1(e^{-1(10)}) = - e^{-1(10)} \,,\quad  B_1(e^{0(10)}) = 0$$
$$B_2(e^{-2}) =0\,,\quad  B_2(e^{-1(10)}) =  - i e^{-1(10)} \,,\quad  B_2(e^{0(10)}) =  - 2 i e^{0(10)}\,.$$
By   Table 1, 
$B_1 =  Ê\ad(E^0_1)|_{\gm}$,  $B_2 = Ê\ad(E^0_2)|_{\gm}$ and the  first claim follows.  The isomorphism between $\wt H^0$ and the subgroup 
$H^0$ of $H \subset \SO_{3,2}^o$ follows by similar  computations that determine explicitly the  elements of  $\GL^{\mathrm{gr}}_0(\gm, J)$  satisfying \eqref{str1}.   
 \end{pf}
 \subsection{The flag of distributions of $P^0$}\hfill\par
 \label{filtrationP^1}
On  $P^0$, there is a  natural flag of distributions 
$T P^0 =  \cD^0_{-2} \supset  \cD^0_{-1} \supset  \cD^0_{(0|-1)} \supset  \cD^0_{(0|0)} \supset \{0\}$,
defined by 
\beq \label{filtration1} \cD^0_{-1} = (\pi^0_*)^{-1}(\cD)\,,\, \cD^0_{(0|-1)} =  (\pi^0_*)^{-1}(\cE)\,,\,  \cD^0_{(0|0)} = (\pi^0_*)^{-1}(\{0\}) = T^{\Vert} P^0\,.\eeq
These distributions  determine the following  filtrations on each $T_y P^0$
\beq \cF\,:\, \underset{ = V_{-2}}{T_y P^0} \supset \underset{ = V_{-1}}{\cD^0_{-1}|_y }\supset \underset{= V_0}{\cD^0_{(0|-1)}|_y}\supset \{0\}\,,\eeq
\beq \cF_*\,:\, \underset{ = V_{-2}}{T_y P^0} \supset \underset{ = V_{-1}}{\cD^0_{-1}|_y }\supset \underset{= V_{0}}{\cD^0_{(0|-1)}|_y}\supsetneq  \underset{= V_{(0|0)}}{\cD^0_{(0|0)}|_y}\supset \{0\}\,,\eeq
so that   $\ggr(T_y P^0, \cF_*)$ is isomorphic as graded vector space to
$$\gm + \gh^0 = \gm^{-2} + \gm^{-1} + \gm^0 + \gh^0 \subset \gg =  \so_{3,2}\,.$$
Any endomorphism  $J_y: \cD^0_{-1}|_y \longrightarrow \cD^0_{-1}|_y$,  which projects onto the complex structure $J_x: \cD|_{x} \longrightarrow  \cD|_{x}$, $x = \pi^0(y)$, 
is a partial complex structure of  $ T_y P^0$ and  two such partial complex structures are equivalent. Due to this,    we may  arbitrarily fix one  such  map $J_y$ at any $y$  and 
  consider the linear frames
\beq \label{adapted1}Êu_\sharp: \gm^{-2} + \gm^{-1} + \gm^0 + \gh^0\longrightarrow V = T_y P^0\eeq
{\it that are adapted to the filtration and partial complex structures of $T_y P^0$}. {\it  This property   does not depend on the choice of $J_y$}. \par
\smallskip
Recall  that a   frame $u_\sharp: \gm \longrightarrow T_y P^0$  is   completely determined  by the corresponding  basis  of  $T_y P^0$
$$f^{-2} = u_\sharp(e^{-2}) \,,\quad f^{-1}_j = u_\sharp(e^{-1}_j)\,,\quad f^{0}_j = u_\sharp(e^{0}_j)\,,\quad F^0_j = u_\sharp(E^0_j)\,.$$
Its  dual  coframe is  indicated by $(f^{-2*}, \ldots , F^{0*}_2)$.\par
\section{The first prolongation   of the Tanaka structure $P^0$}
\setcounter{equation}{0}
\subsection{Adapted frames of $P^0$}\hfill\par
In the next definition, we denote by  $\wh E^0_i$, $i = 1,2$,   the fundamental vector fields determined by   the  right actions of   
$E^0_i  \in Lie(\wt H^0)=\gh^0 $ on $P^0$.
\begin{definition} Let $y \in P^0$ be a point   over  $x = \pi^0(y) \in M$. A linear frame $u_\sharp = (f^{-2}, f^{-1}_i, f^0_j, F^0_k)$, adapted to the filtration and partial complex structure of $T_y P^0$, Êis called {\it 
adapted frame of $P^0$} if 
\begin{itemize} 
\item[i)] $F^0_1 = \wh E^0_1|_y$ and $F^0_2 = \wh E^0_2|_y$;
\item[ii)]  the projected frame $\underline u_\sharp = (\underline f^{-2}, \underline f^{-1}_i , \underline f^{0}_j)  =$  $(\pi^0_*(f^{-2}),  \pi^0_*(f^{-1}_i),  \pi^0_*(f^0_j))$ of $T_x M$ belongs to  the equivalence class $y$, i.e. $y = [\underline u_\sharp]$.
\end{itemize}
The collection $\cF r_{1*}(P^0)$ of  such frames is called {\it bundle of adapted frames of $P^0$} and we denote by  $\pi_{1*}: \cF r_{1*}(P^0)\longrightarrow P^0$   the natural projection. 
\end{definition}
\smallskip
For any adapted frame $u_\sharp: \gm + \gh^0 \longrightarrow T_y P^0$, let     $u: \gm + \gh^0 \longrightarrow \ggr(T_y P^0, \cF)$ be the corresponding isomorphism of graded vector spaces. By  Remark \ref{P^0interpretation},  all    frames  $u_\sharp \in \cF r_{1*}(P^0)|_y$ have the same associated  isomorphism $u$, so that, by Lemma \ref{lem4.2} and the remarks  in \S \ref{section4.3}, we have:\par
\begin{lemma}Ê The triple $(\cF r_{1*}(P^0), P^0, \pi_{1*})$ is a principal bundle of frames over $P^0$, with structure group
$$G^{1*}_{\sharp}Ê= GL_{1*}(\gm + \gh^0, J)\,.$$ 
For any $y \in P^0$, the fiber $\cF r_{1*}(P^0)|_y$ is in natural one-to-one correspondence with the collection of all  s.h.s. $(K^{-2}, K^{-1}, K^{(0|-1)})$ of 
$V = T_y P^0$. 
\end{lemma}
Later on, we will constantly  identify an  {\it adapted frame}  at  $y \in P^0$ with the  corresponding {\it adapted s.h.s.} $K = (K^{-2}, K^{-1}, K^{(0|-1)})$ of 
$V = T_y P^0$.   Moreover,  given an adapted  s.h.s. $K \subset T_y P^0$,  associated with the   linear frame $\wh K = u_\sharp$,  and an element  $X\in  \gm + \gh^0$,  we denote  $$X_{K}  = \wh K(X)   \in T_y P^0\,.$$
 Moreover, if $\wh K'  = \wh K \cdot A$ for some  $A \in GL_{1*}(\gm+\gh^0,  J)$, we write 
\beq\label{tino} X_{K'} = A(X_K)\,,\quad \text{where we denote by}\quad A(X_{K}) = 
\wh {K} \circ A \circ {\wh {K}}^{-1}(X_K)\,.\eeq
\par
\subsection{$\alpha$-torsion, $\beta$-torsion and $c$-torsion}\hfill\par
\label{fundamentaltorsion}
Consider   a smooth field of adapted  s.h.s.  on  a  neighbourhood $\cU$ of  $y \in P^0$
$$\cK : \cU \subset P^0\longrightarrow  \{\,\text{adapted s.h.s. in}\,T_{y'} P^0\,,\,y' \in \cU\,\} \simeq \left.\cF r_{1*}(P^0)\right|_{\cU}\,,$$
We call
 {\it  torsion of  $\cK $ at $y$\/}  the
bilinear map 
\beq \tau_{\cK , y}  \in   \Hom(\Lambda^2\gm, \gm+  \gh^0) \,,\qquad \tau_{\cK , y}(X,Y) = \wh{\cK _y}^{-1}\left(\left.\left[ X_{\cK },
Y_{\cK }\right]\right|_y\right)\,.\label{7.4}\eeq
Notice  that   the  space   $\Tor(\gm) = \Hom(\L^2 \gm, \gg)$,  
into which any torsion $ \tau_{\cK , y}$ takes values,   is  naturally graded, with   homogeneous  subspaces 
\beq \label{7.5} \Tor^k(\gm) = \{ \t \in \Tor(\gm): \t(\gm^i, \gm^j) \subset \gg^{i+j+k}\,,\,\text{for all}\, i,j = -1, -2\}\,,\eeq
where,  as usual, $\gg^{-2} = \gm^{-2}$, $\gg^{-1} = \gm^{-1}$, $\gg^0 = \gm^0 + \gh^0$,  $\gg^1=\gh^1$,
$\gg^2=\gh^2$ and  $\gg^i = \{0\}$ otherwise.  Let us write  $\tau_{\cK , y} $
as  a sum of homogeneous components
$$\tau_{\cK , y}   = \textstyle \sum_{k} \tau^k_{\cK , y}   \qquad\text{with}\qquad   \tau^k_{\cK , y}  \in \Tor^k(\gm)\,.$$
{\it A priori, the torsion $\t_{\cK , y}$ depends  not only on  $ K  = \cK |_y$, but also on Ê  other parts of the  first order jet of $\cK $ at $y$.\/} Nonetheless, there are      components  of $\t_{\cK , y}$  that depend only on $ K = \cK|_y$. 
Such components are very important, because they allow to impose  conditions, which are preserved by CR diffeomorphisms and  determine   canonical   reductions of $\cF r_{1*}(P^0)$. 
In the next definition,  we  name  a few  components that we later show  to  enjoy such crucial property.   \par
\begin{definition} Let  $\cK $ be   a local field of adapted s.h.s. on a neighbourhood $\cU \subset P^0$ of $y$. We call
\begin{itemize}
\item[--] {\it $\a$-torsion  at $y$}  the  map in $ \Hom(\gm^{-1} \times \gm^0, \gm^{-1})$
\beq 
 \a_{\cK , y} = \t^0_{\cK , y}|_{\gm^{-1} \times \gm^0}\,; \eeq
 \item[--]  {\it $\b$-torsion at $y$}  the  map in $\Hom  (\gm^{-1(10)} \times \gm^{0(10)}, \gm^{0(01)})$
\beq 
 \b_{ \cK , y} = (\t^1_{\cK , y}|_{\gm^{-1(10)} \times \gm^{0(10)}} )_{\gm^{0(01)}}\,, \eeq
 where  $(\cdot)_{\gm^{0(01)}}$ denotes the component in the antiholomorphic subspace of $\gm^{0\bC} = \gm^{0(10)} + \gm^{0(01)}$; 
 \item[--] {\it $c$-torsion at $y$} the   map   in $\Hom(\L^2(\gm^{-2} + \gm^{-1}), \gm + \gh^0)$
 \beq c_{\cK , y} = \left.\t_{\cK , y}\right|_{\L^2( \gm^{-2}Ê+ \gm^{-1})}\,.\eeq
\end{itemize}
We denote by $c^k_{\cK , y}$ the homogeneous components   $c^k_{\cK , y} = \left.\t^k_{\cK , y}\right|_{\L^2( \gm^{-2}Ê+ \gm^{-1})}$.
\end{definition}
\subsection{Strongly adapted frames of $P^0$ and the first prolongation $P^1$}\hfill\par
\subsubsection{A preliminary step: the reduction  $\wt{\cF r}_1(P^0)$ of $ \cF r_{1*}(P^0)$}\hfill\par
\label{secondstep}
As before,  $\cK $ is   a  field  of adapted s.h.s. on  a neighbourhood  of $y \in P^0 $.
\begin{lem}  \label{primolemma}Ê The $\alpha$-torsion $\a_{\cK , y}$ depends only on  $K = \cK |_y$ in $T_y P^0$ and 
can be considered as a tensor $\a_K$,   naturally  associated with  $\wh K \in \cF r_{1*}(P^0)$. \par
The collection $\cF r_1(P^0) \subset \cF r_{1*}(P^0)$ of  adapted frames $\wh K $ such that 
\beq \label{condition1}
\left(\a_{K }(e^{-1(10)}, e^{0(10)})\right)_{\gm^{-1(10)}} = 0 = \left(\a_{K }(e^{-1(01)}, e^{0(10)})\right)_{\gm^{-1(01)}}\,,\eeq
 is a  principal subbundle  with structure group 
$\GL_{1}(\gm + \gh^0, J)$. 
\end{lem}
\begin{pf} Let $\cK $,  $\cK '$ be two fields of adapted s.h.s. with $\cK |_y   = \cK '|_y  = K $.  Then there exists a map 
$$A = I + B: \cU \longrightarrow \GL_{1*}(\gm + \gh^0, J)\,,\quad B_{y'} \in \ggl_{1*}(\gm + \gh^0,  J)\,,\, \,y' \in \cU\,,$$
with $B_y = 0$ and  such that $\wh \cK '_{y'} = \wh \cK _{y'} \circ A_{y'}$ for any $y'$ of a neighbourhood of $y$. 
It follows that, for any $X, Y \in \gm + \gh^0$,
\beq\label{6.8} \left.[X_{\cK '}, Y_{\cK '}]\right|_y = \left.[X_{\cK }, Y_{\cK }]\right|_y + 
((X_{\cK }\cdot B|_y)(Y))_{K} - ((Y_{\cK } \cdot B|_y)(X))_{K}\,,\eeq
where we denoted by $X_{\cK }\cdot B|_y$, $Y_{\cK }\cdot B|_y$ the linear maps
$$X_{\cK }\cdot B|_y = \left.\frac{d B(\g_s)}{d s}\right|_{s= 0} \,,\qquad Y_{\cK }\cdot B|_y = \left.\frac{d B(\h_s)}{d s}\right|_{s= 0}\,,$$
for some curves $\g_s$, $\h_s$  with $\g_0 = \h_0 = y$ and $\dot \g_0 = X_{\cK }|_y$, $\dot \h_0 = Y_{\cK }|_s$. Notice that 
$X_{\cK }\cdot B|_y$, $Y_{\cK }\cdot B|_y \in \ggl_{1*}(\gm + \gh^0, J)$. So, if $X \in \gm^{-1}$, $Y \in \gm^0$, 
$$\a_{\cK ', y}(X, Y) = \left(\wh \cK '_y{}^{-1}( \left.[X_{\cK '}, Y_{\cK '}]\right|_y)\right)_{\gm^{-1}} \!\!\!\!\!\!\!\!=  \left(\wh \cK _y^{-1}\left.[X_{\cK }, Y_{\cK }]\right|_y\right)_{\gm^{-1}} \!\!\!\!\!\!\!\!= \a_{\cK , y}(X, Y) $$
proving  the first claim. 
For the second claim, we need  to show that: 
\begin{itemize}
\item[a)] for any $y \in P^0$, there exists an adapted s.h.s. $K $ in $T_y P^0$ for which \eqref{condition1} holds; 
\item[b)]Êtwo such adapted s.h.s. are related one to the other by an element in $\GL_1(\gm + \gh^0, J)$. 
\end{itemize}
 Consider a fixed s.h.s. $K _o$ in $T_y P^0$ and a field of s.h.s. $\cK _o$  on a neighbourhood  of  $y$ with $\cK _o|_y = K _o$. Let also $K $ a second s.h.s. in $T_y P^0$ and a field of s.h.s. $\cK $ with $\cK |_y = K $. As before, we have that $\wh \cK |_{y'} = \wh \cK _o|_{y'} \circ A_{y'}$, $y' \in \cU$,  for some map $A = I + B$ with values in $\GL_{1*}(\gm + \gh^0, J)$. By  definitions  of the distributions $\cD^0_{(0|-1)}  \subset  \cD^0_{(-1)}$ and of the adapted frames,  for any $X \in \gm^{-1}$, $Y \in \gm^0$, the vector  $\a_K (X, Y)$ is equal to 
\beq \label{6.9} \left( \!\wh K ^{-1}\left(\left.[X_{\cK }, Y_{\cK }]\right|_y\right) \!\right)_{\gm^{-1}}\!\!\!\!\!\!\!\!\!=  \left(\!A^{-1} \circ \wh K _o^{-1}\left([X_{\cK _o} + B(X_{\cK _o}), Y_{\cK _o} + B(Y_{\cK _o})]|_y\right)\!\right)_{\gm^{-1}} \!\!\!\!\!\!\!\!=  $$
$$ =   \a_{K _o}(X, Y) + \left((I + B)^{-1}\left(\wh K _o^{-1}\left([X_{\cK _o}, B(Y_{\cK _o}) ]|_y \!\!\! \mod \cD^0_{(0|-1)}\right)\right)\right)_{\gm^{-1}}= $$
$$ =  \a_{K _o}(X, Y) +  \left(\wh K ^{-1}_o\left(\left.[X_{\cK _o}, B(Y_{\cK _o}) ]\right|_y\right)\right)_{\gm^{-1}}\,.\eeq
Since $B$ takes values in $\ggl_{1*}(\gm + \gh^0, J)$, one has that, {\it modulo terms of higher grade},   the values of $B_y$ 
are  of the form
\beq \label{Bdefinition} B_y(e^{-2})  = \l e^{-1(10)} +  \overline \l e^{-1(01)}\,,\quad  B_y(e^{-1(10)})  = \m e^{0(10)} \,\,,$$
$$ B_y(e^{0(10)})  = \nu E^{0(10)} + \nu' E^{0(01)}\,, \qquad \text{for some}\, \l, \m, \nu,  \nu'\in \bC\,.\eeq
Therefore, by \eqref{6.9} we have that
\beq \label{6.10} \a_{K }(e^{-1(10)}, e^{0(10)} ) =  \a_{K _o}(e^{-1(10)}, e^{0(10)} ) + $$
$$ + \left(\wh K _o^{-1}\left(\nu  \left.\left[e^{-1(10)}_{\cK _o}, E^{0(10)}_{\cK _o}\right]\right|_y   +  \nu'  \left.\left[e^{-1(10)}_{\cK _o}, E^{0(01)}_{\cK _o}\right]\right|_y \right)\right)_{\gm^{-1}} = $$
$$ = \a_{K _o}(e^{-1(10)}, e^{0(10)} ) + \nu  e^{-1(10)} \,,\eeq
\smallskip
\beq\label{6.11}Ê \a_{K }(e^{-1(01)}, e^{0(10)} ) = \a_{K _o}(e^{-1(01)}, e^{0(10)} ) +$$
$$ +  \left(\wh K _o^{-1}\left(\nu \left. \left[e^{-1(01)}_{\cK _o}, E^{0(10)}_{\cK _o}\right]\right|_y   +  \nu'  \left.\left[e^{-1(01)}_{\cK _o}, E^{0(01)}_{\cK _o}\right]\right|_y \right)\right)_{\gm^{-1}} = $$
$$ = \a_{K _o}(e^{-1(01)}, e^{0(10)} ) + \nu'  e^{-1(01)} \,,\eeq
where we used the fact that   $E^{0(10)}_{\cK _o}$ and $E^{0(01)}_{\cK _o}$ are fundamental vector fields. 
From \eqref{6.10} and \eqref{6.11}, one can  directly  see that  there always exist  $\nu $,  $\nu'$  such that $K $  satisfies \eqref{condition1}  and that two given s.h.s. $K $, $K '$ satisfy  \eqref{condition1}  if and only if their corresponding adapted frames are related by a transformation $A = I + B$, in which $B$ acts on $e^{0(10)}$  as in \eqref{Bdefinition} with $\nu, \nu' = 0$. Since this is equivalent to say that $B \in \ggl_1(\gm + \gh^0, J)$, we get that  $\cF r_1(P^0)$ is a $\GL_1(\gm + \gh^0,  J)$-reduction. 
\end{pf}
From now on, {\it we  consider only adapted s.h.s. and   fields of adapted s.h.s., whose corresponding adapted frames are in the reduction 
$\cF r_1(P^0)$}.
\begin{lem}\label{secondolemma} If we restrict to $\cF r_1(P^0)$, for any field $\cK $ of adapted s.h.s. around  $y \in P^0$, the  $\b$-torsion  $\b_{\cK , y} $ depends only on $K  = \cK |_y$ and can be considered as a tensor $\b_K $,  naturally  associated with  $\wh K \in \left.\cF r_1(P^0)\right|_y$. \par
The collection 
$\wt{\cF r}_1(P^0) \subset \cF r_1(P^0)$ of   frames     with $\b_K  = 0$  is a  subbundle  with a structure group $G^1_\sharp$, which contains $ \GL_{2}(\gm+\gh^0,  J)$ as normal subgroup and such that $L^1 \= G^1_\sharp/\GL_{2}(\gm+\gh^0,  J)$ is the set of  
 equivalence classes 
\beq L^1 = \{\,I + B\!\!\!\!\mod \GL_2(\gm + \gh^0, J) \,:\phantom{aaaaaaaaaaaaaaaaaaaaaaa}$$
$$\phantom{aaaa} \,  B\in \ggl^{\text{\rm gr}}_1(\gm + \gh^0,  J)\,\,\text{\rm such that}\, \,(B(e^{-1(10)}))_{\gm^{0(10)}}Ê= \m e^{0(10)}\,,$$
$$\phantom{aaaa}\, (B(e^{-1(10)}))_{\gh^{0\bC}}Ê= \nu E^{0(10)} +  (\nu - \overline \m) E^{0(01)} \,\text{\rm for some}\, \nu, \m \in \bC\,\}\,.\label{7.6bis} \eeq 
\end{lem}
\begin{pf} The first claim is proved as  in Lemma \ref{primolemma}. In fact, let $\cK $, $\cK '$ be two  fields of s.h.s. in $\cF r_1(P^0)$, with $\cK |_y = \cK '|_y = K $,  and denote by $A = I + B$ a $\GL_1(\gm + \gh^0, J)$-valued map such that $\wh \cK ' = \wh \cK  \circ A$. Using the same notation of \eqref{6.8}, for any $X \in \gm^{-1(10)}$, $Y \in\gm^{0(01)}$, we have  that 
$ X_{\cK }\cdot B|_y$, $Y_{\cK }\cdot B|_y$ are in  $ \ggl_1(\gm + \gh^0, J)$ and    
$$\left(Y_{\cK }\cdot B|_y(X)\right)_K  \in \widehat {K }(\gm^{0(10)} + (\gh^0)^{\bC})\quad, \quad \left(X_{\cK }\cdot B|_y(Y)\right)_K  = 0\,.$$
From this and \eqref{6.8},  it follows that $(\t^1_{\cK ', y}(X, Y))_{\gm^{0(01)}} = (\t^1_{\cK , y}(X, Y))_{\gm^{0(01)}}$, so   that   $\b_{\cK ', y} = \b_{\cK , y}$. \par
Also the second claim is proved as in   Lemma \ref{primolemma}. Consider a fixed  s.h.s. $K _o$ with $\wh K _o  \in \cF r_1(P^0)|_y$ and let  $\cK _o$ be a field of s.h.s., with  associated frames  in $\cF r_1(P^0)$,  such that  $\cK _o|_y = K _o$.  Take  also  a second s.h.s. $K $ in $\cF r_1(P^0)|_y$ and a field of s.h.s. $\cK $ with $\cK |_y = K $, with corresponding frames   in $\cF r_1(P^0)$. Finally, let  $A = I + B$ be a $\GL_1(\gm + \gh^0, J)$-valued map such that $\wh \cK = \wh \cK _o \circ A$.\par
\smallskip
Since $\pi^0_*(\cD^0_{(0|-1)}) = \cE$ is an integrable complex distribution  of complex dimension one, there is no loss of generality if we assume that Êthe (locally defined)   fields of s.h.s. $\cK $,  
$ \cK '$  are  such  that 
\beq \label{6.14} \pi^0_*([e^{0(10)}_{\cK _o}, e^{0(01)}_{\cK _o}]) = 0 =  \pi^0_*([e^{0(10)}_{\cK }, e^{0(01)}_{\cK }]) \eeq
or, equivalently,  that the vector fields  $[e^{0(10)}_{\cK _o}, e^{0(01)}_{\cK _o}]$, $[e^{0(10)}_{\cK }, e^{0(01)}_{\cK }]$  take values in  $ \cD^0_{(0|0)} = T^{\Vert} P^0$. 
Since  these frames  satisfy  \eqref{condition1}, we also  have that 
\beq\label{6.15} [e^{0(10)}_{\cK _o}, e^{-1(01)}_{\cK _o}] = \r e^{-1(10)}_{ \cK _o} \mod  \cD^0_{(0|-1)}\eeq
for some suitable complex function $\r$. Moreover, by   \eqref{adapted0}, 
$$[e^{-1(10)}_{\cK _o}, e^{-1(01)}_{\cK _o}] = \frac{i}{2} e^{-2}_{\cK _o}Ê\mod   \cD^0_{-1}\,.$$
$$[[e^{0(10)}_{ \cK _o}, e^{-1(01)}_{\cK _o}] , e_{\cK _o}^{-1(01)}] =  - \frac{i}{2}  e^{-2}_{\cK _o} \mod   \cD^0_{-1}\,.$$
This and  \eqref{6.15} imply that  $\rho \equiv -1$ and  that 
\beq \label{eccoqua}Ê[e^{0(10)}_{\cK _o}, e^{-1(01)}_{\cK _o}] = - e^{-1(10)}_{\cK _o}  = [e^{0(10)}, e^{-1(01)}]_{\cK _o} \mod   \cD^0_{(0|-1)}\,.\eeq
These arguments  hold for any field    of s.h.s. in $\cF r_1(P^0)$ and imply  that \eqref{eccoqua}  holds for  $\cK $ as well.
From  \eqref{6.8}, \eqref{6.14} and \eqref{eccoqua} one gets  
\beq  \label{6.17}
 \b_{K }(e^{-1(10)}, e^{0(10)}) = \left(\wh \cK ^{-1}\left(\left.[e^{-1(10)}_{\cK }, e^{0(10)}_{\cK }]\right|_y\right)\right)_{\gm^{0(01)}} = $$
 $$ = \left((I + B_y)^{-1} \circ \wh \cK _o^{-1}\left(\left.[e^{-1(10)}_{\cK }, e^{0(10)}_{\cK }]\right|_y\right)\right)_{\gm^{0(01)}} = $$
$$ =   \b_{K _o}(e^{-1(10)}, e^{0(10)})  - \left(\!B_y([e^{-1(10)}, e^{0(10)}])\right)_{\!\gm^{0(01)}}\!\!\!+ \left(\![B_y(e^{-1(10) }), e^{0(10)}]\right)_{\!\gm^{0(01)}}\eeq
 Since $B$ takes values in $\ggl_{1}(\gm + \gh^0, J)$,  the images   under the linear map $B_y$ are vectors  of the form  (modulo terms of higher grade)
 \beq \label{defB} B_y(e^{-2}) = \l e^{-1(10)} + \overline \l e^{-1(01)}\,,\,\,B_y(e^{-1(10)})Ê= \m e^{0(10)} + \nu E^{0(10)} + \nu' E^{0(01)} \,,$$
$$  B_y(e^{-1(01)})Ê= \overline \m e^{0(01)} + \overline \nu E^{0(01)} + \overline \nu' E^{0(10)} \,, \,B_y(e^{0(10)}) =   B(e^{0(01)}) = 0\eeq
for some  $\l$, $\m$, $\nu$, $\nu'$ $ \in \bC$. So, by  Table 1,   \eqref{6.17} is equivalent to 
\beq \label{stellettabis} 
 \b_{K }(e^{-1(10)}, e^{0(10)}) = \b_{K _o}(e^{-1(10)}, e^{0(10)}) + (- \overline \m + \nu - \nu') e^{0(01)}\,. \eeq
This shows  that if $\nu' = \nu - \overline \m  - \b_{K _o}(e^{-1(10)}, e^{0(10)}) $,  then $\b_K = 0$, so that  $\wt{\cF r}_1(P^0)|_y \neq \emptyset$. From \eqref{stellettabis}, it follows also   that 
$\wt{\cF r}_1(P^0)$  is a  reduction  of $\cF r_1(P^0)$ with structure group $G^1_\sharp$. 
 \end{pf}
\medskip
\subsubsection{The strongly adapted frames of $P^0$}\hfill\par
\label{section6.3.2}
From now on,  we limit ourselves   to  the  bundle $\pi_{\wt 1}: \wt{\cF r}_1(P^0)\longrightarrow P^0$ and {\it the  s.h.s. $K $'s or fields of adapted s.h.s. $\cK $ are  assumed to correspond to   frames in $\wt{\cF r}_1(P^0)$}. Such frames  are called  {\it nicely adapted}. \par
\begin{lem} \label{acca}ÊFor any field $\cK $ of nicely adapted s.h.s. around  $y \in P^0$, the $k$-th components $\t^k_{\cK , y}$ and $c^k_{\cK , y}$ of   torsion  and $c$-torsion  are such that:\par
\begin{itemize}
\item[(i)] $\t^k_{\cK , y} = 0$ for any $k < 0$;
\item[(ii)] $\t^0_{\cK , y}(X,Y) = [X,Y]$  for any  $X\in \gm^{-2} + \gm^{-1}$ and $Y \in \gm$; 
\item[(iii)]  $c^1_{\cK , y} $ depends only on the s.h.s. $K  = \cK |_y$ and it can be considered as a tensor  $c^1_K $,  associated with $\wh K \in \left.\wt{\cF r}_1(P^0)\right|_y$. \par
\end{itemize}
\end{lem}
\begin{pf} For what concerns (i) and (ii), the only statement that does not follow directly from \eqref{adapted0} and the definitions is  claim (ii)  for the case   $X \in \gm^{-2}Ê+ \gm^{-1}$ and $Y \in \gm^0$.  Assume that  $X \in \gm^{-1} $,  $Y \in\gm^0$. Then  \eqref{condition1} and  \eqref{eccoqua} imply   that
$$\left.[X_{ \cK }, Y_{\cK }]\right|_y = \left.([X, Y]_{ \cK })\right|_y \!\!\! \mod  \cD^0_{(0|-1)}|_y\,,$$
from which it follows   $\t^0_{\cK , y}(X, Y) = [X, Y]$. On the other hand, when  $X = e^{-2}$ and $Y = e^{0(10)}$, we have that  
$e^{-2}_{\cK } = - 2i [e^{-1(10)}_{ \cK }, e^{-1(01)}_{ \cK }]  \!\!\! \mod  \cD^0_{(0|-1)}$, so that,  by \eqref{condition1} and  \eqref{eccoqua},  
$$[e^{-2}_{\cK }, e^{0(10)}_{\cK }] =  - 2i[ [e^{-1(10)}_{\cK }, e^{-1(01)}_{\cK }],  e^{0(10)}_{ \cK }]  \!\!\! \mod  \cD^0_{-1} \, = 0 \!\!\!\! \mod  \cD^0_{-1}\,. $$
This shows  that $\t^0_{\cK , y} (e^{-2}, e^{0(10)}) = 0 = [e^{-2}, e^{0(10)}]$ and  (ii) follows. \par
 Now,  as in the proof of  Lemmas \ref{primolemma}  and \ref{secondolemma},  consider a fixed nicely adapted s.h.s. $K _o$  and a field $\cK _o$ of nicely adapted s.h.s. around $y$ with $\cK _o|_y = K _o$.  Take  also  a second nicely adapted s.h.s. $K $  and a field of nicely adapted  s.h.s. $\cK $ with $\cK |_y = K $ and denote by $A = I + B$ the $G^1_\sharp$-valued map such that $\wh \cK = \wh \cK _o \circ A$.  Consider the expression \eqref{6.8} for the Lie bracket $\left.[X_{\cK '}, Y_{\cK '}]\right|_y$ in the case  $X \in \gm^i$, $Y\in \gm^j$, $i, j  \in \{-1,-2\}$. 
 Since $B$ takes values in $Lie(\GL^1_\sharp$) and $Lie(\GL^1_\sharp)/ \ggl_2(\gm+\gh^0, J) = \gl^1$ with $\gl^1 = Lie(L^1)$,  it follows that  
 $$X_{\cK }\cdot B|_y\!\!\!\mod \ggl_2(\gm+\gh^0, J)\,,\qquad Y_{\cK }\cdot B|_y\!\!\!\mod \ggl_2(\gm+\gh^0, J)$$
  are in $\gl^1 $, while  the last two
terms of \eqref{6.8} take value in 
$\left.\cD^0_{j+1}\right|_y$ and $\left.\cD^0_{i+1}\right|_y$, 
 respectively. Since $i+1, j+1 > i + j +1$, it follows that 
$c^1_{\cK ', y}(X,Y) = \left(\wh {K }^{-1}\left(\left.[X_{\cK }, Y_{ \cK }]\right|_y\right)\right)_{\gm^{i+j+1}} = c^1_{\cK ,y}(X,Y)$. 
\end{pf}
\medskip
Recall now that, for  a graded Lie algebra  $\gn$ and a graded $\gn$-module $W$, the
space of skew-symmetric multi-linear maps 
$$C^\ell_k(\gn, W) = \{\,c \in \Hom(\Lambda^\ell\gn, W)\,:\,c(\gn^{i_1}\wedge \dots \wedge\gn^{i_\ell})
\subset W^{i_1 + \dots + i_\ell + k}\,\}$$
 is called  {\it homogeneous space of $\ell$-cochains  of degree $k$}. Its {\it differential\/}  is the
 coboundary
operator $\partial: C^\ell_k(\gn, W) \longrightarrow   C^{\ell+1}_k(\gn, W)$   
$$\partial c(X_1\wedge \dots \wedge X_{\ell+1}) =
\sum_s (-1)^{s+1} X_s \cdot  c(X_1 \wedge \ldots \hat{\underset{s}{\phantom{s}}} \ldots \wedge X_{\ell+1}) + \phantom{aaaaaaaa}$$
$$ \phantom{aaaaaaaaaaaaaaa} +
\sum_{s<t} (-1)^{s+t} c([X_s,X_t] \wedge X_1 \wedge \dots \hat{\underset{s}{\phantom{s}}} 
\ldots \hat{\underset{t}{\phantom{s}}}\ldots \wedge X_{\ell+1})$$
and the corresponding cohomology spaces are denoted by 
$$H^\ell_k(\gn, W) = \frac{\operatorname{Ker} \partial|_{C^\ell_k(\gn, W)}}{\partial(
C^{\ell-1}_k(\gn, W))}\,.$$
Let    $\gm_- = \gm^{-2} + \gm^{-1}$ and $\gh_+ = \sum_{i > 0} \gh^i$.  Consider 
    $\gg = \so_{3,2} =  \gm^{-2}Ê+ \gm^{-1}Ê+ (\gm^0 + \gh^0)+ \gh^1 + \gh^2$ as a graded
$\gm_-$-module and notice that
$$\ggl^{\text{gr}}_k(\gm_- + (\gm^0  + \gh^0) +  \sum_{i = 1}^{k-1} \gh^i) \simeq C^1_k(\gm_-, \gg)   
\,,\qquad  \Tor^k(\gm) \simeq C^2_k(\gm_-,\gg) \,,$$
where 
$\ggl^{\text{gr}}_k(\gm_-  + \sum_{i = 0}^{k-1} \gg^i)  = 
\{B \!\mod \ggl_{k+1}(\gg), \,B \in \ggl_{k}(\gg) \}$.
Hence,  the differential  $\partial$ gives a map  from 
$\ggl^{\text{gr}}_k(\gm_-  + \sum_{i = 0}^{k-1} \gg^i)$ into $ \Tor^k(\gm)$.
Recall also that, via the Killing form  of $\gg$, each space  $C^\ell_k(\gm_-, \gg)$
can be identified with $C^\ell_k(\gh_+, \gg^*)$ and  the opposite of the map $\partial: C^\ell_k(\gh_+, \gg^*) \longrightarrow C^{\ell+1}_k(\gh_+, \gg^*)$ can be identified with  a linear map 
$$\partial^*: C^\ell_k(\gm_-, \gg) \longrightarrow C^{\ell-1}_k(\gm_-, \gg)\,,$$
called  {\it codifferential\/}.
By Kostant's theory (\cite{Ko}),  for any $\ell, k \geq 0$, 
we have  the  $\ad(\gg^0)$-invariant direct
sum decomposition
\beq \label{Kostant} C^\ell_k(\gm_-, \gg) = \partial C^{\ell-1}_k(\gm_-, \gg) \oplus H^\ell_k(\gm_-,\gg)
\oplus \partial^* C^{\ell + 1}_k(\gm_-, \gg)\,,\eeq
from which it follows   the $(\gh_+  + \gg^0)$-invariance of the  spaces 
$$\ker\partial^*|_{C^{\ell}(\gm_-, \gg)}Ê=  \sum_k \left(H^\ell_k(\gm_-, \gg) + \partial^* C^{\ell +1}_k(\gm_-, \gg)\right)$$
(we consider $\gh_+ + \gg^0$ acting on $\gm_- \simeq \gg/(\gh_+ + \gg^0)$ with the adjoint action).\par
\medskip
Consider now the abelian Lie algebra $\gl^1 = Lie(L^1)$ and the  complementary subspace of  $\partial \gl^1$ in $\Tor^1(\gm)$, defined as follows. First of all, 
notice that the restriction $\ad(E^0_2)|_{\gg^i} $  on each   graded space $\gg^i = \gm^i$ or $\gh^i$  is either trivial or equal to a multiple of $J|_{\gg^i}$. This means  that  
$K =  e^{\bR \ad(E_0^2)|_{\gg^i}} $ is trivial or isomorphic to $S^1$ and we may consider a  $K$-invariant Euclidean inner product $<, >$
on  $\partial C^1_1(\gm_-, \gg)$. On the other hand, from definitions, one can directly check that $\gl^1$ is $\ad(E^0_2)$-invariant and, consequently,  the orthogonal complement $(\partial \gl^1)^\perp$  in $ \partial C^1_1(\gm_-, \gg)$ is  $\ad(E^0_2)$-invariant as well. The space $(\partial \gl^1)^\perp$   is also  $\ad(E^0_1)$-invariant, because the   action of $E^0_1$ on $ \partial C^1_1(\gm_-, \gg) $ is equal to   the identity: in fact,    $E^0_1$ is a grading element and 
$\partial C^1_1(\gm_-, \gg) $ is a homogeneous space of  grade $+1$.
So, $(\partial \gl^1)^\perp$  is  $\gh^0$-invariant and complementary to $\partial \gl^1$ in  $\partial C^1_1(\gm_-, \gg)$, while  $(\partial \gl^1)^\perp + \ker\partial^*|_{C^2_1(\gm_-, \gg)}$ is  $\ad(\gh^0)$-invariant and  complementary  to $\partial \gl^1$ in  $\Tor^1(\gm)
\simeq C^2_1(\gm_-,\gg)$. This observation allows the following
\begin{definition}  We call  {\it strongly adapted frame of $T_y P^0$} a nicely adapted frame $\wh K  \in \wt{\cF r}_1(P^0)|_y$, with  $c$-torsion such that 
\beq c^1_K  \in (\partial \gl^1)^\perp + \ker\partial^*|_{C^2_1(\gm_-, \gg)}\,.\eeq 
\end{definition}
\medskip
\begin{lem} \label{lemma6.8} 
The collection  $P^1_\sharp \subset \wt{\cF r}_1(P^0)$  of strongly adapted frames  is a reduction  with a structure 
group  $\wt H^1_\sharp \subset G^1_\sharp$ with  $\GL_{2}(\gm  + \gh^0) \subset \wt H^1_\sharp$ as normal subgroup and 
such that 
$\wt H^1 \= \wt H^1_\sharp/\GL_{2}(\gm  + \gh^0) $ is the set of  equivalence classes 
 \beq\label{6.22} \wt H^1 = \{\, I + B\!\!\!\! \mod \GL_2(\gm + \gh^0)\,: \,  B \in \gl^1 \subset  \ggl^{\text{\rm gr}}_1(\gm + \gh^0,J),\,  \partial B = 0\, \}.\eeq
\end{lem}
\begin{pf} Let  $K _o$ be  a  nicely adapted s.h.s. in $T_y P^0$ and $K $ some other nicely adapted s.h.s. in the same tangent space. As usual, we  consider 
two fields of nicely adapted s.h.s. $\cK _o$,  $\cK $ with $\cK _o|_y = K _o$, $\cK |_y = K $ and a local map  $A = I + B: \cU \longrightarrow  G^1_\sharp$  such that $\wh \cK  = \wh \cK _o \circ A$.  Using  definitions,  Lemma \ref{acca} (i),  (ii) and standard arguments, for  $X \in \gm^{i}$, $Y \in \gm^j$, $i,j \in \{-2, -1\}$, 
\beq \label{6.22bis} c_K ^1(X, Y)  = \left((I + B)^{-1}\left(\wh K ^{-1}_o(\left.[X_{\cK }, Y_{\cK }]\right|_y)\right)\right)_{\gm^{i + j + 1}} = $$
$$ = c_{K _o}^1(X, Y) - \left(B\left(\t^0_{\cK _o, y}(X, Y)\right)\right)_{\gm^{i + j + 1}}  +
\left(\t^0_{\cK _o, y}(B(X), Y)\right)_{\gm^{i + j + 1}}   +$$
$$ +  \left(\t^0_{\cK _o, y}(X, B(Y))\right)_{\gm^{i + j + 1}} =   c^1_{K _o}(X, Y) - \partial \wt B(X, Y)\,,\eeq
where  $\wt B$ is the map in $C^1_k(\gm_-, \gg) \!\simeq\! \ggl^{\text{gr}}_1(\gm + \gh^0)$ such that $I + \wt B = I + B\! \mod \GL_2(\gm + \gh^0)$.
Hence, if we denote by 
$  c^1_{K } = \left(c^1_{K }\right)_{\partial \gl^1} + \left(c^1_{K }\right)_{ (\partial \gl^1)^\perp + \ker\partial^*|_{C^2_1(\gm_-, \gg)}}$  
the  decomposition of $c^1_{K } $  into  a sum of elements in $\partial \gl^1$ and $ (\partial \gl^1)^\perp + \ker\partial^*|_{C^2_1(\gm_-, \gg)}$, we  see that  there always exists  $B$ such that 
$ \left(c^1_{K }\right)_{\partial \gl^1} = 0$, proving that   the fiber  of $P^1_\sharp$ over $y$ is not empty. The equality \eqref{6.22bis}  shows also that  nicely adapted frames $\wh K $, $\wh K '$ are  both in $P^1_\sharp$ if and only if 
$\wh K '= \wh K  \circ A$ for some  $A = I + B$, with $\partial \wt B = 0$. This implies that the structure group $\wt H^1_\sharp$ of $P^1_\sharp$  includes $\GL_2(\gm + \gh^0, J)$ as normal subgroup and is such that $\wt H^1 = \wt H^1_\sharp/ \GL_{2}(\gm + \gh^0, J)$ is as in \eqref{6.22}. Since $\GL_{2}(\gm + \gh^0, J) = \GL_{2}(\gm + \gh^0)$, the claim follows. 
 \end{pf}

\begin{lem} \label{lemma6.9} 
The  abelian Lie algebra $Lie(\wt H^1)$ has dimension $2$ and
$$Lie(\wt H^1) = \{\,\wt X \in \ggl_1^{\text{\rm gr}}(\gm, J)\,:\,\wt X =  \ad(X)|_{\gm}\!\!\!\! \mod \gh\,, \,X \in \gh^1\,\}\,.$$
In particular, $Lie(\wt H^1) $ is isomorphic to $\gh^1$ as vector space and it is naturally  endowed with the complex structure  $J|_{\gh^1}$ defined in \eqref{cstrgg}.
\end{lem}
\begin{pf} Recall that, if we identify the elements of   $\gl^1 = Lie(L^1)$ with linear maps $B \in C_1^1(\gm_-, \gg) \simeq \ggl_1^{\text{gr}}(\gm + \gh^0)$, a linear map $B$ is in $\wt H^1$  if and only if  it  satisfies  \eqref{defB}  for some  $\nu' = \nu - \overline \m$. On the other hand,  condition $\partial B = 0$ is equivalent  to   
$$\left(B([e^{-2}, e^{-1(10)}])\right)_{\gm^{-2}} \!\!\!\!= \left([B(e^{-2}), e^{-1(10)}] + [e^{-2}, B(e^{-1(10)})] \right)_{ \gm^{-2}}\,,$$
$$\left(B([  e^{-1(10)}, e^{-1(01)}])\right)_{\gm^{-1}} \!\!\!\!= \left([B(e^{-1(10)}), e^{-1(01)}] + [e^{-1(10)}, B(e^{-1(01)})]\right)_{\gm^{-1}}\!\!\!\! \,.$$
Using Table 1, one can check that these conditions correspond to require
$$  \nu = \frac{i}{2} \overline \l\,, \qquad \m = - \frac{i}{2} \l\,,\qquad \nu' = 0\,, $$
from which it follows   that $Lie(\wt H^1)$  is generated by
   the (equivalence classes of the) maps 
   $B_1$, $B_2$  corresponding to  $\l = 1$ and $\l = i$, respectively, i.e.,   
\beq\label{4.24} B_1(e^{-2}) = e^{-1(10)} + e^{-1(01)}  = \ad(- E^1_2)(e^{-2})\,,$$
$$B_1(e^{-1(10)})Ê= - \frac{i}{2} e^{0(10)}  +\frac{i}{2}  E^{0(10)}  = \ad(-E^1_2)(e^{-1(10)})\,,$$
$$B_1(e^{0(10)}) =  0 =  \ad(-E^1_2)(e^{0(10)})\!\!\!\mod \gg^1\,,\eeq 
\beq \label{4.28}ÊB_2(e^{-2}) = i (e^{-1(10)} - e^{-1(01)})  = \ad( E^1_1)(e^{-2})\,,$$
$$ B_2(e^{-1(10)})Ê= \frac{1}{2} e^{0(10)} +  \frac{1}{2} E^{0(10)}  = \ad(E^1_1)(e^{-1(10)}) \,,$$
$$B_2(e^{0(10)}) =  0 =  \ad(E^1_1)(e^{0(10)})\!\!\!\mod \gg^1 \,,\eeq
and this concludes the proof.
\end{pf}
\medskip
The quotient bundle  
$$\pi^1: P^1 =  P^1_\sharp/GL_{2}(\gm+\gh^0)\longrightarrow P^0$$
with  $\pi^1$ induced by the natural projection $\pi^1_\sharp: P^1_\sharp\longrightarrow P^0$, is called {\it first prolongation of the Tanaka structure $\pi^0: P^0 \longrightarrow M$}. It is a principal
bundle over $P^0$,  but it also a principal bundle over $M$. In fact, 
\begin{lem} \label{lemma6.10} There exists a natural  right action of a semidirect product
 $\wt H^0 \ltimes \wt H^1$ on the  bundle $\pi^0 \circ \pi^1: P^1 \longrightarrow M$, which makes  it a 
 principal bundle over $M$, canonically associated with the  girdled CR structure $(\cD, J)$.\par
 The  Lie group $\wt H^0 \ltimes \wt H^1$ is isomorphic to  $H/H^2$, where we denote by  $H^2 $ 
 the connected subgroup of $H = \Aut_{x_o}(M_o) \subset \SO_{3,2}^o$ with subalgebra $Lie(H^2) = \gh^2$.
\end{lem}
\begin{pf} For any $h \in H^0$, we denote by $f_h: P^0 \longrightarrow P^0$ the  diffeomorphism determined  by the right action on $P^0$, and by 
$\wh f_h: \cF r(P^0) \longrightarrow \cF r(P^0)$ the associated  diffeomorphism on the linear frame  bundle $\cF r(P^0)$  of $P^0$,   defined by 
 $$\wh f_h(u_\sharp) = f_h{}_* \circ u_\sharp \circ \Ad_{h}: \gm + \gh^0 \longrightarrow T_{f_h(y)} P^0$$
   for  any $u_\sharp \in \cF r(P^0)|_y$.
 Using the   definition of the adapted frames, one can check that $\wh f_h$ maps  $\cF r_{1*}(P^0)$ into itself. Moreover, using the fact that $\wh f_h$ preserves also the  distributions  \eqref{filtration1}  and the partial complex structures of the tangent spaces,  one can check  that  $\wh f_h(\cF r_1(P^0)) \subset \cF r_1(P^0)$ and  
 $\wh f_h(\wt{\cF r}_1(P^0)) \subset \wt{\cF r}_1(P^0)$.\par
  Now,  from Lemma \ref{lemma5.4} and the fact that $(\partial \gl^1)^\perp + \ker\partial^*|_{C^{\ell}(\gm_-, \gg)}$ is $\ad(\gh^0)$-invariant, it follows that $\wh f_h$ maps strongly adapted frames into strongly adapted frames, inducing an automorphism $\wt f_h: P^1 \longrightarrow P^1$. This shows the existence of a right action of $\wt H^0$  on $P^1$ and, consequently,  of a  right action of a semidirect product 
  $\wt H^0 \ltimes \wt H^1$, which  acts transitively  and freely on the fibers of $\pi^0 \circ \pi^1: P^1 \longrightarrow M$. 
The last claim can be checked for instance computing the product rule of the Lie group 
$\wt H^0 \ltimes \wt H^1 \simeq \bC^* \ltimes \bC$ and comparing it with the group structure of $H/H^2$, 
determined by  products of  matrices in $H\subset\SO_{3,2}^o$ 
(it is convenient to  use the represention in  \S \ref{section3.2}). 
The group structure of   $\wt H^0 \ltimes \wt H^1$ can be determined using  Lemmas \ref{lemma5.4} and \ref{lemma6.8}.  
\end{pf}
In analogy with  Remark \ref{P^0interpretation},  if
 $u_\sharp, u'_\sharp: \gm + \gh^0 \longrightarrow T_y P^0$ are two  strongly adapted  frames in  the same   equivalence class $[u_\sharp] \in P^1|_y$, they determine the same  graded  isomorphism  $u=u': \gm + \gh^0 \longrightarrow \ggr(T_y P^0, \cF)$. \par
 \medskip
We conclude observing that, in  analogy with  \S \ref{filtrationP^1}, there is a natural flag of distributions on $P^1$,  given by  
$$ \cD_{-1}^1 = (\pi^1_*)^{-1} (\cD_{-1}^0),\quad
\cD_{(0|-1)}^1 = (\pi^1_*)^{-1}( \cD_{(0|-1)}^0),\quad
\cD_{(0|0)}^1 = (\pi^1_*)^{-1}( \cD_{(0|0)}^0)$$
and $\cD^1_1 = (\pi^1_*)^{-1}(0) = T^{\Vert} P^{1}$. These distributions and the CR structure of $(M, \cD, J)$ induce filtrations of type \eqref{filtration}, \eqref{filtrationstar} on each tangent space
\beq \cF\,:\, \underset{ = V_{-2}}{T_y P^1} \supset \underset{ = V_{-1}}{\cD^1_{-1}|_y }\supset \underset{= V_0}{\cD^1_{(0|-1)}|_y}\supset \underset{= V_1}{\cD^1_{1}|_y}\supset \{0\}\,,\eeq
\beq \cF_*\,:\, \underset{ = V_{-2}}{T_y P^1} \supset \underset{ = V_{-1}}{\cD^1_{-1}|_y }\supset \underset{= V_{0}}{\cD^1_{(0|-1)}|_y}\supsetneq  \underset{= V_{(0|0)}}{\cD^1_{(0|0)}|_y}\supset  \underset{= V_{1}}{\cD^1_{1}|_y}\supset \{0\}\eeq
endowed with a (unique up to equivalences)  partial complex structure $J$. By construction, the graded vector spaces $\ggr(T_yP^1,  \cF_*)$ 
and $\gm + \gh^0 + \gh^1$ are isomorphic.
\par
\begin{rem} \label{remark6.11}
 Let $\wt E^0_i$, $\wt E^1_i$, $i,j = 1,2$,   be the four 
fundamental   vector fields of $\pi^1: P^1\longrightarrow P^0$, 
determined by
the right (infinitesimal) actions  of  $E^0_i$, $E^1_j$ $ \in  
Lie(\wt H^0 \ltimes \wt H^1)$  (which is naturally isomorphic,  as vector space,  
to $\gh^0 + \gh^1$) and set 
$\wt E^{\ell(10)} = \frac{1}{2}(\wt E^\ell_1 - i \wt E^\ell_2)$ for 
$\ell = 0,1$.  For any   collection of (real and complex) vector 
fields
$(\wt e^{-2} , \wt e^{-1(10)}, \wt e^{-1(01)}, \wt e^{0(10)}, \wt 
e^{0(01)} )$ on some open set $\cU \subset P^1$, such that, for any $z \in 
\cU$, the projected vectors
$e^{-2}_z =  \pi^1_*(\wt e^{-2}_z)$, $e^{-1(10)} = \pi^1_*(\wt 
e^{-1(10)_z})$, etc.,  are associated with   a strongly adapted frame in  $ T_{\pi^1(z)} P^0$, one has the following identities(which will be used in the next section):
\beq \label{6.28} [\wt E^{1(10)}, \wt e^{0(10)}] = 0\ ,\qquad [\wt 
E^{1(10)}, \wt e^{0(01)}] = - \wt E^{1(01)}\ .\eeq
 They can be directly inferred from  the Jacobi identity for the Lie brackets between $ \wt e^{0(01)}$ and $[\wt E^{1(10)}, \wt e^{0(10)}]$ and between $ \wt e^{0(01)}$ and 
  $[\wt 
E^{1(10)}, \wt e^{0(01)}]$. 
\end{rem}
\section{The second prolongation   of the Tanaka structure $P^0$}
\setcounter{equation}{0}
\subsection{Adapted frames  of $P^1$}\hfill\par
\begin{definition} Let $z \in P^1|_y$ be a point  over  $y = \pi^1(z) \in P^0$. A linear frame $u_\sharp: \gm + \gh^0 + \gh^1 \longrightarrow T_z P^1$, adapted to the filtration and partial complex structure of $T_z P^1$, Êis called {\it 
adapted frame of $P^1$} if 
\begin{itemize} 
\item[i)] the restriction $ u_\sharp|_{\gh^0 + \gh^1}: \gh^0+\gh^1 \longrightarrow \left.\cD^1_{(0|0)}\right|_z $ coincides with the isomorphism given by the right action of 
$ Lie(\wt H^0 \ltimes \wt H^1) = Lie(H/H^2)$ ($\simeq  \gh^0 + \gh^1$  as vector space) on  $ P^1$;
\item[ii)]  the  projected  frame $\underline u_\sharp = \pi^1_* \circ  u_\sharp|_{\gm + \gh^0}: \gm + \gh^0 \longrightarrow T_y P^0$ is in the equivalence class $z = [\underline u_\sharp] \in P^1$.
\end{itemize} \par
The collection of  such frames is called {\it bundle of adapted frames of $P^1$} and is denoted by $\cF r_{2*}(P^1)$. We denote by  $\pi_{2*}: \cF r_{2*}(P^1) \longrightarrow P^1$  the natural projection.  
\end{definition}
\smallskip
By remarks in \S \ref{section4.3}  and  \S \ref{section6.3.2},  any linear frame of $\cF r_{2*}(P^1)|_z$ is  uniquely determined by the corresponding s.h.s. $(K ^{-2}, K ^{-1}, K ^{(0|-1)}, K ^{(0|0)})$ of 
$V = T_z P^1$.  From this and usual  arguments (see also Remark \ref{triviality}),   it follows that:
\begin{lemma}Ê The triple $(\cF r_{2*}(P^1), P^1, \pi_{2*})$ is a principal bundle  over $P^1$, with structure group $GL_{2*}(\gm + \gh^0 + \gh^1)$. 
\end{lemma} 
In  analogy with \S \ref{fundamentaltorsion},  for  a given smooth field $\cK $ of adapted s.h.s. in the tangent spaces of a neighbourhood  of $z \in P^1$, we may consider the {\it torsion of $\cK $ at $z$}
\beq \tau_{\cK , z}  \in   \Hom(\Lambda^2\gm, \gm+  \gh^0 + \gh^1) \,,\qquad \tau_{\cK , z}(X,Y) = \wh{\cK _z}^{-1}\left(\left.\left[ X_{\cK },
Y_{\cK }\right]\right|_z\right)\,, \label{7.1}\eeq
 the associated  {\it $c$-torsion}  $c_{\cK , z} = \t_{\cK , z}|_{\L^2(\gm^{-2} + \gm^{-1})}$ and the graded components $\t^k_{\cK , z}$ and $c^k_{\cK , z}$. Besides this, we need to consider the following 
\begin{definition} Given  a local field $\cK $ of adapted s.h.s. on a neighbourhood  of $z$, we call
{\it $\g$-torsion  at $z$} the restriction  
$$\g_{\cK , z} = \left.\tau^1_{\cK , z}\right|_{\gm^{-1} \times \gm^0} \in \Hom  (\gm^{-1} \times \gm^0, \gm^0 + \gh^0).$$ 
\end{definition}
\begin{lem}  \label{lemma7.4} 
The $\gamma$-torsion $\g_{\cK , z}$ depends only on  $K  = \cK |_z$ in $T_z P^1$ and 
can be considered as a tensor $\g_K $   associated with   $\wh K $. 
\par
The subset $\cF r_2(P^1) \subset \cF r_{2*}(P^1)$ of  adapted frames $\wh K $ such that 
\beq \label{condition1P2}
\left(\g_{K }(e^{-1(10)}, e^{0(10)})\right)_{\gm^{0(10)}} = 0 = \left(\g_{K }(e^{-1(01)}, e^{0(10)})\right)_{\gm^{0(01)}}  \eeq
 is a  reduction  with structure group 
$\GL_{2}(\gm + \gh^0 + \gh^1)$. 
\end{lem}
\begin{pf} The proof of the first claim is the perfect analogue of the argument used for Lemma \ref{primolemma}. Also the second claim is proved in a very similar way.  In fact, consider a fixed s.h.s. $K _o$ in $T_z P^1$ and a field of s.h.s. $\cK _o$  on a neighbourhood $\cU$ of  $z$ with $\cK _o|_z = K _o$. Any other field of s.h.s. $\cK $ is such that $\wh \cK |_{z'} = \wh \cK_o |_{z'} \circ A_{z'}$, $z' \in \cU$,  for some map $A = I + B$ with values in $\GL_{2*}(\gm + \gh^0 + \gh^1)$.  By usual arguments,  we find that    $\g_K (X, Y)$,  with $K  = \cK _z$, $X \in \gm^{-1}$, $Y \in \gm^0$, is equal to
\beq \label{7.4bis}\g_{K }(X, Y)  =   \g_{K _o}(X, Y) +  \left([X, B_z(Y) ]\right)_{\gm^{0}}\,.\eeq
On the other hand, modulo terms of higher  grades, the  map $B_z$ is 
such that
$$B_z(e^{-2})  = \l e^{0(10)} +  \overline \l e^{0(01)}\,,\quad  B_z(e^{-1(10)})  = \m E^{1(10)} + \mu' E^{1(01)}\,,$$
$$ B_z(e^{0(10)})  = \nu E^{1(10)} + \nu' E^{1(01)}\qquad \text{for some}\,\,\l, \m, \m',  \nu,  \nu'\in \bC\,,$$
from which it follows that 
$$\g_{K }(e^{-1(10)}, e^{0(10)} ) =  \g_{K _o}(e^{-1(10)}, e^{0(10)} ) - \frac{1}{2}Ê\nu  e^{0(10)} - \frac{1}{2} \nu' E^{0(10)}\,,$$
$$\g_{K }(e^{-1(01)}, e^{0(10)} ) =  \g_{K _o}(e^{-1(01)}, e^{0(10)} ) - \frac{1}{2}Ê\nu  E^{0(01)} - \frac{1}{2} \nu' e^{0(01)}\,. $$
Hence, there always exist    $\nu $,  $\nu'$  such that $\cK $ satisfies \eqref{condition1P2}, so that   $\cF r_1(P^1)|_z$ is not empty. 
The same expressions show that  $\cF r_2(P^1)$ is a reduction with  structure group $\GL_{2}(\gm + \gh^0 + \gh^1)$. 
\end{pf}
\smallskip
\subsection{Strongly adapted frames of $P^1$  and the  second prolongation}\hfill\par
The proof of next lemma is basically the same of Lemma \ref{acca} (iii). \par
\begin{lem} \label{accaP2}ÊLet  $\cK $ be a field of s.h.s. on a neighbourhood of $z$, with associated adapted frames  in the  reduction $\cF r_2(P^1)$ defined in Lemma \ref{lemma7.4}.  
Then   $c^2_{\cK ,z}$ depends only on $K  = \cK |_z$ and it can be considered as a tensor $c^2_K $, associated with $\wh K  \in \left.\cF r_2(P^1)\right|_z$. 
\end{lem}

Now, as it is pointed out in \S \ref{section6.3.2}, we have that 
$$ \ggl_2^{\text{gr}}(\gm + \gh^0 + \gh^1) \simeq C^1_2(\gm_-, \gg)\,,\qquad \Tor^2(\gm) \simeq C^2_2(\gm_-, \gg)\,.$$
So, by \eqref{Kostant}, we have  that   $\ker\partial^*|_{C^2_2(\gm_-, \gg)}$
 is   complementary  to $\partial\ggl_2^{\text{gr}}(\gm + \gh^0 + \gh^1) $ in  $\Tor^2(\gm)
\simeq C^2_2(\gm_-,\gg)$ and it  is invariant under the adjoint actions $\ad_X$, for $X\in \gh^0 + \gh^1$. 
\par
\begin{definition}  We call  {\it strongly adapted frame of $T_z P^1$} any  frame $\wh K  \in  {\cF r}_2(P^1)|_z$, 
whose  $c$-torsion is such that 
\beq Êc^2_K  \in  \ker\partial^*|_{C^2_2(\gm_-, \gg)}\,.\eeq 
\end{definition}
 
\begin{lem} \label{lemma7.6} 
The subset  $P^2_\sharp \subset {\cF r}_2(P^1)$  of strongly adapted frames  is a reduction  with a structure 
group $\wt H^2_\sharp$, which contains $\GL_3(\gm + \gh^0 + \gh^1)$ as normal subgroup and such that
$\wt H^2\= \wt H^2_\sharp/ \GL_3(\gm + \gh^0 + \gh^1)$ is the set of equivalence classes 
 $$\wt H^2 = \{\,I + B\!\!\!\!\mod\GL_3(\gm + \gh^0 + \gh^1)\,:\, B \in \ggl^{\text{\rm gr}}_2(\gm + \gh^0 + \gh^1) \,,\partial B = 0\,\}\,.$$
 \end{lem}
 The proof of  Lemma \ref{lemma7.6}  is  identical to  Lemma \ref{lemma6.8} and we omit it. \par
\begin{lem}\label{kappa} The real Lie algebra $\wt \gh^2 = Lie(\wt H^2)$ is $1$-dimensional and
 generated  by 
 $B \in C^1_2(\gm_-, \gg) \simeq \ggl_2^{\text{\rm gr}}(\gm + \gh^0 + \gh^1)$,  defined by 
\beq \label{B1P2}Ê
B(e^{-2}) =  E^{0(10)} + E^{0(01)} , \   B(e^{-1(10)})Ê= i E^{1(10)},\   B(e^{-1(01)})Ê= -i  E^{1(01)}. \eeq
\end{lem}
\begin{pf}  
A map  $B \in \ggl_2^{\text{\rm gr}}(\gm + \gh^0 + \gh^1)$  is  of the form
$$B(e^{-2}) = \l e^{0(10)} + \overline \l e^{0(01)} + \mu E^{0(10)} + \overline{\mu} E^{0(01)}\,,$$
$$B(e^{-1(10)})Ê= \nu E^{1(10)} +  \nu' E^{1(01)}  \,,\quad  B(e^{-1(01)})Ê= \overline \nu' E^{1(10)} +  \overline \nu E^{1(01)} $$
for some $\l$, $\mu$ $\nu$, $\nu'$  $ \in \bC$. The condition $\partial B = 0$ is equivalent to   
$$\left(\!B([e^{-2}, e^{-1(10)}])\!\right)_{\gm^{-1}} = \left([B(e^{-2}), e^{-1(10)}] + [e^{-2}, B(e^{-1(10)})] \right)_{ \gm^{-1}}\,,$$
$$\left(\!B([  e^{-1(10)}, e^{-1(01)}])\!\right)_{\!\gm^{0} + \gh^0} \!\!\!\!\!\!\!\! =  \left(\![B(e^{-1(10)}), e^{-1(01)}] + [e^{-1(10)}, B(e^{-1(01)})]\!\right)_{\!\gm^{0} + \gh^0} \!\!.$$
A direct computation shows that this holds if and only if
  $$ i \mu =  \nu\,,\quad  i \overline \mu =   \nu\,,\quad  - i  \overline  \l = \nu'  \,\quad     i \overline \l  =  \nu' \,,$$
or, equivalently, if and only if
  $(\l, \m, \nu, \nu') = (0, t,  i t, 0)$ with $ t \in \bR$. The value $t=1$ gives  the generator  defined by \eqref{B1P2}.
\end{pf}
By Lemma \ref{kappa},  $Lie(\wt H{}^2)$ is 1-dimensional and  is generated by the map $B$, described  in \eqref{B1P2}. 
Notice  that $B$ is equal to the linear map
\beq \label{7.10}  B(X) = \ad_{E^2}(X)\!\!\!\! \mod \gh^2\quad \text{for any}\ X \in \gm + \gh^0 + \gh^1 = \gg/\gh^2\,,\eeq
so that the linear map $\imath:  Lie(\wt H{}^2) \longrightarrow \gh^2$, $\imath(B) = E_2$, is a  vector space isomorphism between 
$Lie(\wt H{}^2)$ and
 $\gh^2$.
 \par
We can now define the  {\it second prolongation of the Tanaka structure $P^0$}. 
It is the quotient bundle 
$$\pi^2: P^2 =  P^2_\sharp/GL_3(\gm+\gh^0 + \gh^1)\longrightarrow P^1,$$
where  $\pi^2$ is  the map induced by the natural projection $\pi^2_\sharp: P^2_\sharp\longrightarrow P^1$. It is  a  principal
bundle over $P^1$, but it is also a principal bundle over $M$. In fact, 
\begin{lem} \label{lemma7.11} There exists a natural  right action of $H = \Aut_{x_o}(M_o)$  on the  bundle $\pi = \pi^0 \circ \pi^1\circ \pi^2: P^2 \longrightarrow M$, which makes  $(P^2, M, \pi)$ an $H$-principal bundle,  canonically associated with the  girdled CR structure of $M$.
\end{lem}
\begin{pf} For any $[h] \in H/H^2$ ($\simeq \wt H^0 \ltimes \wt H^1$),  let $f_{[h]}: P^1 \longrightarrow P^1$ be the corresponding diffeomorphism, given by the right action, and by 
$\wh f_{[h]}: \cF r(P^1) \longrightarrow \cF r(P^1)$ the  diffeomorphism of   $\cF r(P^1)$ defined for any $u_\sharp \in \cF r(P^1)|_z$ by  
 $$\wh f_{[h]}(u_\sharp) = f_h{}_* \circ u_\sharp \circ \wt \Ad_{[h]}:  \gm + \gh^0 + \gh^1 = \gg/\gh^2 \longrightarrow T_{f_h(z)} P^1\,,$$
where $\wt \Ad_{[h]}$ denotes the map $\wt \Ad_{[h]}(X) = \Ad_h(X) \!\!\mod \gh^2$,  for any $X \in  \gm + \gh^0 + \gh^1 $.    Using definitions and Remark \ref{remark6.11},  one can check that $\wh f_{[h]}$ maps $\cF r_{2*}(P^1)$ into itself and, since $\wh f_{[h]}$ preserves  the flag of distributions and the partial complex structures of the tangent spaces,  it  is such that $\wh f_{[h]}(\cF r_2(P^1)) \subset \cF r_2(P^1)$. \par
Using Lemmas \ref{lemma5.4}, \ref{lemma6.9} and \ref{lemma6.10} and the fact that $ \ker\partial^*|_{C^{\ell}(\gm_-, \gg)}$ is $\ad(\gh^0 + \gh^1)$-invariant, we obtain that $\wh f_{[h]}$ 
preserves  the bundle of strongly adapted frames $P^2_\sharp$ and induces an automorphism of the $\wt H{}^2$-bundle  
$\pi^2: P^2 \longrightarrow P^1$. 
This shows the existence of a right action of $H/H^2 = \wt H^0\ltimes \wt H^1$ on $P^2$ and one can consider a  
map (not a  group homomorphism)
$$\wt \rho: H/H^2 \times \wt H^2 \simeq (\wt H^0\ltimes \wt H^1) \times \wt H{}^2 \longrightarrow \operatorname{Diff}(P^2)\ ,$$
whose restrictions 
$\wt \rho|_{H/H^2}Ê\longrightarrow  \operatorname{Diff}(P^2)$, $\wt \rho|_{\wt H^2}Ê\longrightarrow  \operatorname{Diff}(P^2)$
 are the group homomorphisms  given by the two right actions of $H/H^2$ and $\wt H^2$  on $P^2$. By construction,  for any $w$, $w' $  in the same fiber of $P^2$ over $x \in M$, there exists a pair $([h], h') \in H/H^2 \times \wt H^2$ such that 
 $\wt \rho(([h], h'))(w) = w'$.
\par
Now, we observe that any  $h \in H \subset \SO_{3,2}^o$ can be uniquely written as 
$$h = h^0 \cdot \exp(X^1) \cdot \exp(X^2)\,,
\quad\text{with}\    h^0 \in H^0,\,  X^1 \in \gh^1,\, X^2 \in \gh^2,$$
where $H^0$ is the connected subgroup of $H$ with $Lie(H^0) = \gh^0$.
A direct way to check this is to use the explicit description of  \S \ref{section3.2} for the matrices in $\gh  \subset \so_{3,2}$. 
We may therefore consider the map
$$\rho: P^2 \times H \longrightarrow P^2, 
\quad
\rho(z, h) = \wh f_{[h^0 \cdot \exp(X^1)]}(z) \cdot \exp(X^2),$$
which one can directly check to be a right action 
that is free and transitive on the fibers of $P^2$ over $M$.  
 \end{pf}
\medskip
 In perfect analogy with $P^1$, there exists  a natural flag of distributions on $P^2$,  given by 
$$ \cD_{-1}^2\= (\pi^2_*)^{-1} (\cD_{-1}^1),\
\cD_{(0|-1)}^2\= (\pi^2_*)^{-1}( \cD_{(0|-1)}^1),\
\cD_{(0|0)}^2\= (\pi^2_*)^{-1} (\cD_{(0|0)}^1),$$
$$\cD_{1}^2\= (\pi^2_*)^{-1} (\cD_1^1),\
\cD^2_2 = (\pi^2_*)^{-1}(0) = T^{\Vert} P^{2}\,.$$
These distributions and the CR structure of $(M, \cD, J)$ determine filtrations of type \eqref{filtration}, \eqref{filtrationstar} 
and a partial complex structure $J$ (determined up to equivalences) on each tangent space $T_w P^2$. All this makes any  graded vector  space $\ggr(T_wP^2,  \cF_*)$  isomorphic to  the graded Lie algebra $ \so_{3,2} = \gm + \gh$.\par
\section{The  Cartan connection of a girdled CR manifold\\and the solution of the equivalence problem}
\setcounter{equation}{0}
Following the same steps  for   the constructions of the first and second prolongations,  we now consider the following
\begin{definition} \label{definition8.1} Let $w \in P^2|_z$ be a point  over  $z = \pi^2(w) \in P^1$. A linear frame $u_\sharp: \gm + \gh \longrightarrow T_w P^2$, adapted to the filtration and partial complex structure of $T_w P^2$, Êis called {\it 
adapted frame of $P^2$} if 
\begin{itemize} 
\item[i)] the restriction $u_\sharp|_{\gh}: \gh \longrightarrow \cD^2_{(0|0)}|_w$
coincides with the isomorphism determined by the right action of $\gh$ on  $ P^2$; 
\item[ii)]  the  projected linear frame $\underline  u_\sharp = \pi^2_* \circ  u_\sharp|_{\gm + \gh^0 + \gh^1}: \gm + \gh^0 + \gh^1 \longrightarrow T_w P^1$ is in the equivalence class $w = [ \underline u_\sharp] \in P^2$.
\end{itemize} \par
The collection of  such frames is called {\it bundle of adapted frames of $P^2$} and is denoted by $\cF r_{3*}(P^2)$. We denote by  $\pi_{3*}: \cF r_{3*}(P^2) \longrightarrow P^2$  the natural projection.  
\end{definition}
\smallskip
By construction,    the  frames in a fiber $\cF r_{3*}(P^2)|_w$ are  completely determined by the corresponding s.h.s. $(K ^{-2}, K ^{-1}, K ^{(0|-1)}, K ^{(0|0)}, K ^1)$ of 
$V = T_w P^2$.  Moreover,  from definitions,    $(\cF r_{3*}(P^2), P^2, \pi_{3*})$ is a principal bundle  over $P^2$ with structure group $\GL_{3*}(\gm + \gh)$. For  any  smooth field $\cK $ of adapted s.h.s.  on a neighbourhood of $w \in P^2$, we may consider the {\it torsion of $\cK $ at $w$}
\beq \tau_{\cK , w}  \in   \Hom(\Lambda^2\gm, \gm+  \gh) \,,\qquad \tau_{\cK , w}(X,Y) = \wh{\cK _w}^{-1}\left(\left.\left[ X_{\cK },
Y_{\cK }\right]\right|_w\right)\,, \label{7.11}\eeq
 the   {\it $c$-torsion}  $c_{\cK , w} = \t_{\cK , w}|_{\L^2(\gm^{-2} + \gm^{-1})}$ and the  following
\begin{definition} Given  a local field $\cK $ of adapted s.h.s. on a neighbourhood $\cU \subset P^2$ of $w$, we call
 {\it $\varepsilon$-torsion} the restriction $\varepsilon_{\cK , w} = \left.\tau^2_{\cK , w}\right|_{\gm^{-2} \times \gm^0}$. 
\end{definition}
\begin{lem}  The $\varepsilon$-torsion $\varepsilon_{\cK , w}$ depends only on  $K  = \cK |_w$ in $T_w P^2$ and 
can be considered as a tensor $\varepsilon_K $,  associated with the   frame $\wh K $. 
\par
The collection  $\cF r_3(P^2) \subset \cF r_{3*}(P^2)$ of  adapted frames  such that 
\beq \label{condition1P3}
\left(\varepsilon_{K }(e^{-2}, e^{0(10)})\right)_{\gh^0} = 0  \eeq
 is a  subbundle  with structure group 
$\GL_3(\gm + \gh)$. 
\end{lem}
\begin{pf} The proof of the first claim is  exactly as in Lemma \ref{primolemma}. For the second claim, as usual,  consider a fixed s.h.s. $K _o$ in $T_z P^1$ and a field of s.h.s. $\cK _o$  on a neighbourhood $\cU$ of  $w$ with $\cK _o|_w = K _o$. Any other field of s.h.s. $\cK $ is such that $\wh \cK |_{w'} = \wh \cK_o |_{w'} \circ A_{w'}$, $w' \in \cU$,  for some map $A = I + B$ with values in $\GL_{3*}(\gm + \gh)$. We have that  $\varepsilon_K (X, Y)$,  with $K  = \cK _w$, $X \in \gm^{-2}$, $Y \in \gm^{0(10)}$, is equal to
\beq \label{7.4ter}\varepsilon_{K }(X, Y))  =   \varepsilon_{K _o}(X, Y) +  \left([X, B_w(Y) ]\right)_{\gm^{0}}\,.\eeq
Modulo a term  of higher  grade, the  map $B_w$ is such that 
$$B_w(e^{-2})  = \l E^{1(10)} +  \overline \l E^{1(01)}\,,\quad  B_w(e^{-1(10)})  = \m E^{2(10)} + \mu' E^{2(01)}\,,$$
$$ B_w(e^{0(10)})  = \nu E^2\,, \qquad \text{for some}\,\,\l, \m, \m',  \nu \in \bC\,.$$
From this,   it follows that 
$$\varepsilon_{K }(e^{-2}, e^{0(10)} ) =  \varepsilon_{K _o}(e^{-2}, e^{0(10)} ) - \nu (E^{0(10)} + E^{0(01)})\,,$$
which can be used to infer that   $\cF r_3(P^2)|_w$ is not empty and that  $\cF r_3(P^2)$ is a reduction with  structure group $\GL_3(\gm + \gh)$. 
\end{pf}
\medskip
The  proof of the following lemma is basically the same  of Lemma \ref{acca} (iii)
and we omit it. \par
\begin{lem}  Let  $\cK $ be a field of s.h.s. on a neighbourhood of $w$, with corresponding   frames in  $\cF r_3(P^2)$.  Then  $c^3_{\cK ,w}$ depends only on $K  = \cK |_w$ and can be considered as a tensor $c^3_K $, associated with $\wh K  \in \left.\cF r_3(P^2)\right|_w$. 
\end{lem}
As we pointed out before, 
$$\ggl_3^{\text{gr}}(\gm + \gh)   \simeq C^1_3(\gm_-, \gg)\,,\qquad \Tor^3(\gm) \simeq C^2_3(\gm_-, \gg)\,,$$
so that  $\ker\partial^*|_{C^2_3(\gm_-, \gg)}$
is complementary  to $\partial\ggl_3^{\text{gr}}(\gm + \gh^0 + \gh^1) $ in  $\Tor^3(\gm)
\simeq C^2_3(\gm_-,\gg)$ and  invariant under $\ad(\gh)$.  
\par
We call    {\it strongly adapted frame of $T_w P^2$} any   adapted frame $\wh K  \in \cF r_{3}(P^2)|_w$  with 
\beq Êc^3_K  \in  \ker\partial^*|_{C^2_3(\gm_-, \gg)}\,.\eeq 
 The proof of the next lemma is essentially the same of  Lemmas \ref{lemma6.8} and \ref{lemma7.11}. \par
\begin{lem} 
The subset 
$P^3_\sharp \subset \cF r_3(P^2)$ of strongly adapted frames is a reduction with structure group 
$H^3_\sharp$ with $ \GL_4(\gm + \gh)$ as normal subgroup and such that $\wt H^3 \= H^3_\sharp/\GL_4(\gm + \gh)$
is the set of equivalence classes of linear maps 
 $$\wt H^3 = \{\,I + B\!\!\! \mod\GL_4(\gm + \gh) \,:\, B \in \ggl^{\text{\rm gr}}_3(\gm + \gh^0 + \gh^1),\ \partial B = 0\,\}\,.$$
 Moreover, there exists a natural  right action of 
 $H\ltimes \wt H^3$
  on the fiber bundle $\pi = \pi^0 \circ \pi^1\circ \pi^2 \circ \pi^3 : P^3 = P^3_\sharp/\GL_4(\gm + \gh) \longrightarrow M$, which makes  $(P^3, M, \pi)$ a 
 principal bundle,  canonically associated with the  girdled CR structure of $M$.\par
 \end{lem}
 Before  our main theorem, it remains only to prove the next lemma,  which  shows that   $\pi^3: P^3  \longrightarrow P^2$ has actually  trivial fibers and that  $P^3 \simeq P^2$ as  $H$-bundles.
 \par
\begin{lem}\label{kappaP3} \label{lemma8.6} The  Lie algebra $Lie(\wt H^3)$ is trivial and  the  structure group   $\pi^3_\sharp: P^3_\sharp \longrightarrow P^2$ is isomorphic to  
$\GL_4(\gm + \gh)$. In particular,  the principal bundle $(P^3, M, \pi)$  is $H$-equivalent to $(P^2, M, \pi)$. 
\end{lem}
\begin{pf} 
Recall that an element $B \in \ggl_3^{\text{gr}}(\gm + \gh)$  is such that  
$$B(e^{-2}) = \l E^{1(10)} + \overline \l E^{1(01)}\,,\quad B(e^{-1(10)})Ê= \mu E^2  \,,\quad B(e^{-1(01)})Ê=   \overline \mu E^{2} $$
for some $\l$, $\mu$  $ \in \bC$. The condition $\partial B = 0$ is equivalent to 
$$\left(B([e^{-2}, e^{-1(10)}])\right)_{\gm^{0} + \gh^0} = \left([B(e^{-2}), e^{-1(10)}] + [e^{-2}, B(e^{-1(10)})] \right)_{\gm^{0} + \gh^0}\,,$$
which  implies that 
$$0 =  \frac{\l}{2}Êe^{0(10)} + \frac{\overline \l}{2} E^{0(10)}Ê- \mu E^{0(10)}Ê - \mu E^{0(01)}\,,$$
 i.e.,  $\l = 0 = \m$. From this, the claim follows.
\end{pf}
We can now proceed with Êthe proof of our main theorem. ÊThe main argument can be outlined as follows. 
As Êit occurs Êin Êthe standard Tanaka construction of Cartan connections (\cite{Ta3, CS, AS}),
iterating Êthe arguments of Êprevious steps, Êone Êgets two natural sequences Êof Êprincipal Êbundles, namely Êthe bundles of Êstrongly adapted linear frames Ê $\pi^{k+1}_\sharp: P^{k+1} _\sharp \longrightarrow P^{k}$, $k \geq 2$, with structure groups $\GL_{k+2}(\gm + \gh)$, 
and the bundles of equivalence classes of linear frames Ê$\pi^{k+1}: P^{k+1} = P^{k+1}_\sharp/\GL_{k+2}(\gm + \gh) \longrightarrow P^k$, Ê each of them trivially equivalent to the others. Since the bundles of linear frames 
$\pi^{k+1}_\sharp: P^{k+1} _\sharp \longrightarrow P^{k}$ have progressively smaller structure groups, they reduce to $\{e\}$-structures for  Ê$k$'s sufficiently high. In our case, 
this occurs when Ê$k+1 \geq 4$, so that  Ê$\pi^4_\sharp: P^{4}_\sharp\longrightarrow P^3$ is Êin fact the collection of linear frames Êof an absolute parallelism on Ê$P^3 (\simeq P^2)$. This is the parallelism that  corresponds to Êthe ÊCartan connection we are looking for.
 \par
\begin{theo}  There exists a Cartan connection $\o: P^2 \longrightarrow \so_{3,2}$ on the principal bundle $(P^2, M, \pi)$,  which is 
 canonically associated with the girdled CR structure,   i.e., satisfies the following two properties: 
 \begin{itemize}
 \item[a)] for any local CR diffeomorphism $f: \cU \subset M \longrightarrow M$,
 the  naturally associated  lifted map $\wt f: \pi^{-1}(\cU) \subset  P^2 \longrightarrow P^2$  is such that $\wt f^* \o = \o$; 
 \item[b)] if  $F: \cV \subset P^2 \longrightarrow P^2$ is a local diffeomorphism such that $F^* \o = \o$, then  $F = \wt f|_\cV$ for some  lifted map $\wt  f$ of a local CR diffeomorphism $f$ of $M$. 
 \end{itemize}
 \end{theo}
\begin{pf} By usual arguments,  we may consider a bundle $\pi_{4}:  \cF r_{4}(P^3)  \longrightarrow P^3$, given by the frames of $P^3$,  defined  in complete analogy with the bundle  described in  Definition \ref{definition8.1}. For any field $\cK $ of s.h.s.,  associated with frames in $\cF r_{4}(P^3)$,  the usual arguments show that 
the $4$-order component $c^4_{\cK , w}$ of the $c$-torsion of $\cK $ at $w \in P^3 \simeq P^2$ depends just on $K  = \cK |_w$ and we may  consider the reduction of $P^4_\sharp \subset \cF r_{4}(P^3)$ given by the  frames  with  $\partial^* c^4_K   = 0$. As in previous proofs,  for any $w \in P^3 (\simeq P^2)$,  two frames in   $P^4|_w$ are  always related by an endomorphism  $A = I + B$ with $B \in \ggl^{\text{gr}}_4(\gm + \gh)$ such that   $\partial B = 0$.  
A simple check shows that $B $ is trivial.  By  the fact that $\ggl_5(\gm + \gh) = 0$, this means that  any fiber of  $\pi_{\sharp}^4 : P^4_\sharp \longrightarrow P^3$ contains exactly   one  element and  that  there exists a unique  section
\beq \label{sectionsigma} \s: P^3  \longrightarrow  \cF r_{4}(P^3) \qquad \text{satisfying }\,\partial^* c^4_{\s}Ê\equiv 0\eeq
$$(\,\text{i.e. with}\ \s_w \in P^4_\sharp|_w\ \text{for any}\  w \in P^3\,)\,.$$  \par 
\smallskip
Since, by  Lemma  \ref{lemma8.6},  the $H$-bundle $\pi: P^2\longrightarrow M$ can be identified with  the $H$-bundle $\pi: P^3 \longrightarrow M$, we may consider 
 the $\gg$-valued 1-form $\o$ on $P^2(\simeq P^3)$, defined by 
$$\o_w = (\s_w)^{-1}: T_w P^2 \simeq T_w P^3 \longrightarrow \gm + \gh$$
at any $w \in P^2 \simeq P^3$. 
We claim that $\o$ is a Cartan connection. In fact, by definitions,  for any $w \in P^2 \simeq P^3$
the linear map 
$$\left.(\o_w)^{-1}\right|_{\gh} = \left. \s_w\right|_\gh : \gh \longrightarrow T^{\Vert}_w P^3 \simeq T^{\Vert}_w P^2$$
coincides with the natural isomorphism between $\gh = Lie(H)$ and  $T^{\Vert}_w P^2$, determined by the right action of $H$ on the fibers. Moreover, by the same arguments of Lemmas \ref{lemma6.10} and \ref{lemma7.11},  the bundle of linear frames $P^4_\sharp$ is invariant under a right action of $H$ on $P^3 \simeq P^2$, so that, for any $h \in H$ and $w \in P^3 \simeq P^2 $, 
$$ R_{h*} \circ \s_w \circ \Ad_h = \s_{w \cdot h}\qquad \Longleftrightarrow \qquad (R_h^* \o)|_w = \Ad_{h^{-1}} \cdot  \o_w\,, $$
proving that $\o$ is  a Cartan connection   modelled on  $M_o =  \SO_{3,2}^o/H$. \par
\medskip
To check  (a), we recall that, by construction, the bundle $\pi^0: P^0 \longrightarrow M$ and the bundles $\pi^{i+1}: P^{i+1} \longrightarrow P^{i}$, $0 \leq i \leq 3$, are   quotients of   bundles  of linear frames of the underlying  manifolds. Using just the definitions, one can  check that  the differential $f_*$ of a CR  diffeomorphism $f$ of $M$ maps the  frames in $P^0_\sharp$ into frames in $P^0_\sharp$ and that the quotient $P^0 = P^0_\sharp/\GL_1(\gm, J)$ is mapped  into itself. This defines a canonical lift $f_*: P^0 \longrightarrow P^0$ on $P^0$. Similarly,   the differential $(f_*)_*$ induces a lift  $(f_*)_*: P^1 \longrightarrow P^1$, the differential $((f_*)_*)_*$ determines a  lift  on $P^2$ and  so on. In particular, we obtain a natural lift $\wt f: P^3 (\simeq P^2) \longrightarrow P^3$, which preserves the unique section $\s: P^3 \longrightarrow P^4_\sharp$ and, consequently,  the Cartan connection $\o$. \par
For (b), consider the basis $\cB = (e^i_j, E^\ell_m)$ of $\so_{3,2} = \gm + \gh$ introduced in \S \ref{section3.2} and the  vector fields $(\wh e^i_j, \wh E^\ell_m)$ defined by 
$$ \wh e^i_j|_w \= (\o_w)^{-1}(e^i_j)\,,\qquad \wh E^\ell_m|_w \= (\o_w)^{-1}(E^\ell_m)\qquad \text{for any} \,w \in P^2\,.$$
Notice that, by the construction of $\o$,  the  flows of the vector fields  $ \wh E^\ell_m$  are    $\Phi^{\wh E^\ell_m}_t = R_{\exp(t E^\ell_m)}$.
Assume that $F: \cV \subset P^2 \longrightarrow P^2$ is a local diffeomorphism such that $F^* \o = \o$ and hence 
such that 
\beq \label{equivariance1}ÊF_*(\wh e^i_j) = \wh e^i_j \,,\qquad F_*(\wh E^\ell_m) = \wh E^\ell_m\,.\eeq
It follows  that $F \circ R_{\exp (t E^\ell_m)} =  R_{\exp (t E^\ell_m)} \circ F$, $t \in \bR$,  and it induces a 
local diffeomorphism $f : \pi(\cV) \subset M \longrightarrow M$ on $M$.  Moreover,  using the fact that  
$$ \cD_{x} = < \pi_*(\wh e^{-1}_i|_w), \pi_*(\wh e^0_i|_w), i = 1,2> \,,\qquad x \in M\,,\,\,w \in \pi^{-1}(x) \subset P^2\,,$$
 one gets  that $f$ is a local CR diffeomorphism. Finally, if we denote by  $\wt f: P^2|_{\pi(\cV)} \longrightarrow P^2$ the  natural lift of  $f$ on $P^2$,  the map  $F' = \wt f^{-1} \circ F: \cV \longrightarrow P^2$
\begin{itemize}
\item[--] is a local diffeomorphism  mapping  the  fields $\wh e^i_j$ and $\wh E^\ell_m$ into themselves; 
\item[--] induces the  identity map $\operatorname{Id}_{\pi(\cV)}$ on $\pi(\cV) \subset M$. 
\end{itemize}
By  the properties of  Cartan connections (use e.g. normal coordinates  -- see \cite{SS1}),   this  occurs if and only if $F'(w) = w$ for any $w \in \cV$, i.e. $F = \wt f|_\cV$.
\end{pf}
Let $\o$ be  the Cartan connection introduced  in  previous  theorem and  
$$\vartheta^{-2}\,,\,\vartheta^{-1(10)}\,,\, \vartheta^{-1(01)}\,,\,\vartheta^{0(10)}\,,\, \vartheta^{0(01)}\,,\, \o^{0(10)}\,,\, \o^{0(01)}\,, \, \o^{1(10)}\,,\, \o^{1(01)}\,,\,  \o^2$$
 the ($\bR$- and $\bC$-valued) 1-forms of $P^2$ that, for any vector field $X$ of $P^2$, give the components  of  the elements  $\o(X) \in \so_{3,2}$  w.r.t. the  basis  formed  by the matrices 
 $e^{-2}$, $e^{-i(10)}$,  $E^{i(10)}$, $E^2$ and  their complex conjugates. Using  the construction of $\o$, one can check that they satisfy the  structure equations
 \beq d \vartheta^{-2} + \frac{i}{2} \vartheta^{-1(10)} \wedge \vartheta^{-1(01)} - \left(\o^{0(10)} + \o^{0(01)}\right)\wedge \vartheta^{-2}= \Theta^{-2}\,,$$
 $$d \vartheta^{-1(10)} - \vartheta^{0(10)}\wedge \vartheta^{-1(01)} - \o^{0(10)}  \wedge \vartheta^{-1(10)} +  i \o^{1(10)}\wedge \vartheta^{-2}  = \Theta^{-1(10)},$$
  $$d \vartheta^{0(10)} - \left(\o^{0(10)} -\o^{0(01)}\right)\wedge \vartheta^{0(10)} + \frac{1}{2} \o^{1(10)} \wedge \vartheta^{-1(10)} = \Theta^{0(10)}\,,$$
$$d \o^{0(10)} - \vartheta^{0(10)} \wedge \vartheta^{0(01)}+ \frac{1}{2} \o^{1(01)} \wedge \vartheta^{-1(10)} + 
\o^2 \wedge \vartheta^{-2} = \O^{0(10)}\,,$$
$$d \o^{1(10)} - \o^{1(01)}\wedge \vartheta^{0(10)} - \o^{1(10)}\wedge\o^{0(01)}  + i \o^2 \wedge \vartheta^{-1(10)} = \O^{1(10)}\,,$$
$$d \o^2 -  \frac{i}{2} \o^{1(10)} \wedge \o^{1(01)} + \left(\o^{0(10)} +  \o^{0(01)}\right) \wedge  \o^2 = \O^{2}\,,\eeq
where  the 2-forms  $\Theta^\a$ and $\O^a$ are  of the  form  (here $\a, \b, \g$  denote  indices as $-2$, $-1(10)$ etc., and  $a, b, c $  denote indices as $2$, $1(10)$, etc.):
$$\Theta^\a = \sum_{\b, \g} T^\a_{\b |\g} \vartheta^\b \wedge \vartheta^\g \,,\qquad \O^a =  \sum_{\b,\g} R^a_{\b |\g} \vartheta^\b \wedge \vartheta^\g $$
with smooth functions $T^\a_{\b|\g}$ and $R^a_{\b| \g}$, called {\it structure functions},  which  satisfy  constraints, corresponding to the conditions considered  in Lemmas \ref{primolemma}, \ref{secondolemma},  \ref{lemma6.8}, \ref{lemma7.4} and \ref{lemma7.6}. For example,  the condition considered in Lemma \ref{primolemma} implies that  
$$T^{-1(10)}_{-1(10)| 0(10)} = T^{-1(01)}_{-1(01)| 0(10)} = 0\,.$$
Analogous constraints come from the other  conditions: each of them either requires  the vanishing of some structure function  or imposes a    linear relation between some of them.

\end{document}